\documentclass[12pt]{amsart}
\usepackage{amsmath, amssymb,verbatim,appendix,enumerate}
\usepackage{hyperref}
\usepackage[mathscr]{eucal}
\usepackage{amscd}
\usepackage{amsthm}
\usepackage{stmaryrd}
\usepackage{comment}
\usepackage[latin1]{inputenc} 
\usepackage{tikz}
\usetikzlibrary{shapes,arrows}
\usepackage{cite}
\usepackage{url}
\usepackage{float}

\usepackage{showlabels}

\usepackage{setspace,a4wide,color,xcolor,graphicx}

\newtheorem{theorem}{Theorem}[section]
\newtheorem{lemma}[theorem]{Lemma}
\newtheorem{corollary}[theorem]{Corollary}
\newtheorem{proposition}[theorem]{Proposition}

\newtheorem{observation}{Observation}
\newtheorem{claim}[theorem]{Claim}
\newtheorem{fact}[theorem]{Fact}

\numberwithin{equation}{section} 

\theoremstyle{definition}

\newtheorem{definition}[theorem]{Definition}

\def\cH{\mathcal{H}}

\def\VC{\mathrm{VC}}

\def\e{\epsilon}
\def\a{\alpha}

\newcommand{\vbar}{\bar{v}}

\newcommand{\Ebar}{\bar{E}}

\newcommand{\calF}{\mathcal{F}}
\newcommand{\calH}{\mathcal{H}}
\newcommand{\calC}{\mathcal{C}}

\newcommand{\calA}{\mathcal{A}}

\newcommand{\calP}{\mathcal{P}}

\newcommand{\calE}{\mathcal{E}}
\newcommand{\calG}{\mathcal{G}}

\newcommand{\calQ}{\mathcal{Q}}

\newcommand{\calS}{\mathcal{S}}

\newcommand{\calU}{\mathcal{U}}

\newcommand{\SVC}{{\mathrm{SVC}}}

\def\Forb{\operatorname{Forb}}

\newcommand{\reg}{{\mathrm{reg}}}
\newcommand{\sml}{{\mathrm{small}}}

\title[]{On the growth of weak regular partitions in k-uniform hypergraphs}
\date{\today }
\author{H. Sheats}
\author{C. Terry}\thanks{The second author was partially supported by NSF CAREER Award DMS-2115518 and a Sloan Research Fellowship.}

\address{Department of Mathematics, Statistics and Computer Science, University of Illinois Chicago, Chicago, IL, USA}
\email{hshea3@uic.edu}
\email{caterry@uic.edu}

\begin{document}

\begin{abstract}
This paper classifies hereditary properties of $k$-uniform hypergraphs according to bounds in the weak hypergraph regularity lemma first developed by Chung. Specifically, we prove the bounds must grow at a constant, polynomial, exponential, or tower rate. This generalizes work in the 3-uniform case by the second author. 
\end{abstract}

\maketitle

\section{Introduction}
The goal of this paper is to determine the possible asymptotic growth rates of weak regular partitions in $k$-uniform hypergraphs, on the level of hereditary properties.  To ease exposition throughout this introduction, we refer to $k$-uniform hypergraphs as simply \emph{$k$-graphs}. 

As is now well-known, there exist many inequivalent generalizations of graph regularity in the $k$-uniform setting. The topic of this paper is the version referred to as \emph{weak hypergraph regularity}, first developed by Chung and Haviland-Thomason \cite{chung, HT}.\footnote{See \cite{part3, part4, GSW2} for similar problems about a stronger form of hypergraph regularity.}  Roughly speaking, given a $k$-graph $H=(V,E)$ and $X_1,\ldots, X_k\subseteq V$, we say $(X_1,\ldots,X_k)$ is \emph{weakly $\e$-regular with respect to $H$} when the density on $(X_1,\ldots, X_k)$ is approximately the same as the density on any sufficiently large subtuple.  The corresponding regularity lemma, proved by Chung, tells us that for any $k$-graph $H=(V,E)$, there is a partition $\calP$ of $V$ of bounded size, such that the regular $k$-tuples from $\calP$ cover almost all of $V^k$  (see Section \ref{sec:defs} for precise definitions).

\begin{theorem}[Chung \cite{chung}]\label{thm:chung}
For all $\e\in (0,1]$ there exists an integer\footnote{The \emph{tower function} is defined by setting $Tw(1)=1$ and $Tw(x+1)=2^{Tw(x)}$.} $M\leq Tw(O_k(\e^{-k-3}))$ so that any $k$-graph $G$ admits a weak $\e$-regular partition into at most $M$ parts.
\end{theorem}

The bound appearing in Theorem \ref{thm:chung} can be deduced from Chung's proof \cite{chung}, and is a tower type bound similar to that arising in the proof of Szemer\'{e}di's regularity lemma for graphs \cite{szemeredi}.  However, it is now also well known that for certain restricted classes of $k$-graphs, much better bounds can be obtained. The first such results appeared in the semi-algebraic setting  \cite{semi-alg, APPRS, FGLNP}, followed by \cite{FPS, CS, CS-distal, AFP}. This leads naturally to the  question of what types of hypergraphs require what types of bounds. As first observed by Alon-Fox-Zhao \cite{AFZ}, the natural setting in which to make this question precise is that of hereditary properties.

A \emph{hereditary $k$-graph property} is a class of finite $k$-graphs closed under isomorphism and induced sub-$k$-graphs.  There is a rich literature exploring dichotomies for hereditary properties based on the behavior of various growth functions (see e.g. \cite{BBW-speed, BBW-penultimate, BT-hypergraphs, Terry-speed, LT-speeds}).  In this paper, we study growth functions which measure  the size of weak regular partitions.  

\begin{definition}
Given a hereditary $k$-graph property $\calH$, define $M_{\calH}:(0,1]\rightarrow \mathbb{N}$ as follows.  For all $\e\in (0,1]$, let $M_{\calH}(\e)$ be the smallest integer $M\geq 1$ so that for all sufficiently large $H\in \calH$, there is a weak $\e$-regular partition of $H$ into at most $M$ parts.
\end{definition}

In \cite{part2}, the second author posed the problem of classifying all possible asymptotic growth rates of functions of the form $M_{\calH}(\e)$, when $\calH$ is a hereditary property of $3$-uniform hypergraphs.\footnote{Here and throughout the paper, the ``asymptotic growth rate" of $M_{\calH}$ means ``as $\e\rightarrow 0$".}  That paper, in conjunction with the main result of \cite{part1}, followed by the improvement on \cite{part1} proved in \cite{GSW}, gives the following classification.

\begin{theorem}[Classification for $3$-graphs \cite{part2,part1,GSW}]
Let $\calH$ be a hereditary $3$-graph property.  Then one of the following  holds.
\begin{enumerate}
\item (tower) There are $C_1,C_2>0$ so that for all sufficiently small $\e>0$,  
$$
Tw(\e^{-C_1})\leq M_{\calH}(\e)\leq Tw(\e^{-C_2});
$$
\item (exponential) There are $C_1,C_2>0$ such that for all sufficiently small $\e>0$,  
$$
2^{\e^{-C_1}}\leq M_{\calH}(\e)\leq 2^{\e^{-C_2}};
$$
\item (polynomial) There are $C_1,C_2>0$ such that for all sufficiently small $\e>0$,  
$$
\e^{-C_1}\leq M_{\calH}(\e)\leq \e^{-C_2};
$$
\item (constant) There is  $C>0$ so that for all sufficiently small $\e>0$,  
$$
M_{\calH}(\e)=C.
$$
\end{enumerate}
\end{theorem}

The goal of this paper is to prove the analogous theorem for $k$-graphs, for all $k\geq 3$ (the case $k=2$ is dealt with in   \cite{part2}).   

\begin{theorem}[Classification for $k$-graphs]\label{thm:main}
Let $k\geq 3$ be an integer, and let $\calH$ be a hereditary $k$-graph property.  Then one of the following  holds.
\begin{enumerate}
\item (tower) There are $C_1,C_2>0$ so that for all sufficiently small $\e>0$,  
$$
Tw(\e^{-C_1})\leq M_{\calH}(\e)\leq Tw(\e^{-C_2});
$$
\item (exponential) There are $C_1,C_2>0$ such that for all sufficiently small $\e>0$,  
$$
2^{\e^{-C_1}}\leq M_{\calH}(\e)\leq 2^{\e^{-C_2}};
$$
\item (polynomial) There are $C_1,C_2>0$ such that for all sufficiently small $\e>0$,  
$$
\e^{-C_1}\leq M_{\calH}(\e)\leq \e^{-C_2};
$$
\item (constant) There is $C>0$ so that for all sufficiently small $\e>0$,  
$$
M_{\calH}(\e)=C.
$$
\end{enumerate}
\end{theorem}

Moreover, the proof of Theorem \ref{thm:main} will show the growth classes are finitely determined in the following sense. There exist for each $x\in \{1,2,3,4\}$, sequences $(\calF^+_n(x))_{n\in \mathbb{N}}$ and $(\calF^-_n(x))_{n\in \mathbb{N}}$, where each $\calF^+_n(x)$ and $\calF^-_n(x)$ is a finite  (possibly empty) set of $k$-graphs, such that $\calH$ is in class $(x)$ if and only if $\calH\cap \calF^-_n(x)\neq \emptyset$ for all integers  $n\in \mathbb{N}$, and on the other hand, for some $m\in \mathbb{N}$, $\calH\cap \calF^+_m(x)=\emptyset$.  Specifically, this tells us that when $M_{\calH}$ grows at a sub-tower rate, there is a finite set of forbidden $k$-graphs responsible.  

Interestingly, Theorem \ref{thm:main} shows the general growth rates in the $k$-uniform setting for $k\geq 4$ are exactly the same as in the $3$-uniform setting.  In contrast, as was shown in \cite{part2}, the $2$-uniform setting allows for only three growth classes: constant, polynomial, and tower.

The proof of Theorem \ref{thm:main} uses $k$-ary generalizations of the combinatorial arguments from \cite{part2}, the new results of  Gishboliner, Shapira, and Wigderson  mentioned above \cite{GSW}, work by Chernikov-Starchenko \cite{CS} and Fox-Pach-Suk \cite{FPS} on hypergraphs of bounded VC-dimension, work of Chung \cite{chung} and Haviland-Thomason \cite{HT}, and the refinement on Gowers' lower bound construction \cite{gowers} proved by Fox-Lov\'{a}sz \cite{FL}.

\vspace{2mm}

\noindent{\bf AI Usage:} The bulk of this paper originates from 2025 and early 2026, before the authors were utilizing AI.  The decision to include Lemma \ref{lem:deletion} was made later in summer 2026, at which point the authors used ChatGPT 5.6 Sol to generate an initial proof. The output was substantially shorter than the human-generated proof in  the 3-uniform case appearing in \cite{part2}.  This new proof was then completely rewritten by the authors and included in Section 7.  All other proofs in the paper were devised  by the authors. ChatGPT 6 Pro was used for proofreading and formatting the bibliography entries. The authors take full responsibility for the results in this paper. 

\section{Preliminaries}\label{sec:defs}
This section covers several necessary preliminaries.  In Subsection \ref{ss:not} we cover basic notation. Subsection \ref{ss:reg} introduces weak hypergraph regularity and homogeneity.  Subsection \ref{ss:HPs} sets up definitions related to hereditary properties. Subsection \ref{ss:blowups} defines blowups, and uses them to state a characterization of when one hereditary $k$-graph property is close to another.  

\subsection{Basic Notation}\label{ss:not}
 Given an integer $\ell\geq 1$, let $[\ell]=\{1,\ldots, \ell\}$.  For a set $X$ and an integer $\ell\geq 0$, define ${X\choose \ell}=\{Y\subseteq X: |Y|=\ell\}$. For an integer $\ell \geq 2$ and sets $V_1,\ldots, V_\ell$ (not necessarily distinct), let
$$
K_\ell[V_1,\ldots, V_\ell]=\left\{v_1\ldots v_\ell\in {V_1\cup \ldots \cup V_\ell\choose \ell}: v_1\in V_1,\ldots, v_\ell\in V_\ell\right\}.
$$

For an integer $k\geq 2$, a \emph{$k$-uniform hypergraph}  is a pair $(V,E)$ where $V$ is a set of vertices and $E\subseteq {V\choose k}$ is a set of edges.  All vertex sets appearing in this paper are assumed to be finite.  For the rest of the paper, we will refer to $k$-uniform hypergraphs as \emph{$k$-graphs}. Given a $k$-graph $G$, we let $V(G)$ denote its vertex set and $E(G)$ its edge set.  We will write $v_1\ldots v_k$ to denote the $k$-element set $\{v_1,\ldots, v_k\}$.  As usual, a \emph{graph} is just a $2$-graph.  A \emph{clique} is a $k$-graph $(V,E)$ where $E={V\choose k}$.  An \emph{anti-clique} is a $k$-graph $(V,E)$ where $E=\emptyset$.

Given a $k$-graph $G=(V,E)$ and $A\subseteq V$, we let $G[A]$ denote the $k$-graph $(A, E\cap {A\choose k})$.  An \emph{induced sub-$k$-graph} of $G$ is any $k$-graph of the form $G[A]$ for $A\subseteq V$. We say $A$ forms a clique or anti-clique in $G$ if $G[A]$ is a clique or anti-clique, respectively.

A $k$-graph $G=(V,E)$ is \emph{$\ell$-partite} if there exists a partition $V=V_1\cup \ldots \cup V_{\ell}$ such that for all $e\in E$ and $1\leq i\leq \ell$, $|e\cap V_i|\leq 1$.     We write $G=(V_1\cup \ldots\cup V_{\ell},E)$ to indicate that $G$ is $\ell$-partite with compatible partition $V_1,\ldots, V_{\ell}$.

When $G$ is a graph, we say it is \emph{bipartite} if it is $2$-partite.  We say $G$ is \emph{cobipartite} if there exist disjoint cliques $V_1,V_2$ so that $V(G)=V_1\cup V_2$. Similarly, we say $G$ is a \emph{split graph} if there is a clique $V_1$ in $G$ and an anti-clique $V_2$ in $G$ such that $V(G)=V_1\cup V_2$.

Suppose $G=(V,E)$ is a $k$-graph. We will frequently work with the set of ordered tuples arising from the unordered edges of $G$.  For this we use the following notation:  
$$
\overline{E}:=\{(v_1,\ldots, v_k)\in V^k: v_1\ldots v_k\in E\}.
$$
  Given nonempty subsets $V_1,\ldots, V_k$ (not necessarily distinct), set 
$$
d_G(V_1,\ldots, V_k)=\frac{|\overline{E}\cap (V_1\times \ldots\times V_k)|}{|V_1|\cdots |V_k|}.
$$
Given a $k$-graph $G=(V,E)$, we let $G^1=G$ and let $G^0$ be the $k$-graph $G^0=(V, {V\choose k}\setminus E)$. Similarly, let $E^1=E$ and $E^0={V\choose k}\setminus E$.  Given $1\leq i<k$ and $e\in {V(G)\choose i}$, the \emph{neighborhood of $e$ in $G$} is 
$$
N_G(e)=\left\{e'\in {V(G)\choose k-i}: e\cup e'\in E(G)\right\}.
$$
When $e=\{v\}$ is a single vertex, we denote this by $N_G(v)$.  We will also use the following more general ``neighborhood" notation for subsets of Cartesian products.  Suppose $X_1,\ldots, X_t$  are sets, $F\subseteq \prod_{j=1}^tX_j$, $1\leq i<t$, and $x\in X_1\times \cdots \times X_i$. Then 
$$
N_F(x)=\{y\in X_{i+1}\times \cdots \times X_t: (x,y)\in F\}.
$$ 

We will frequently use the following standard lemma (see Lemma 2.16 in \cite{part2} for a proof). 

\begin{lemma}\label{lem:averaging}
Let $a,b,\e\in (0,1)$ be such that $ab=\e$. Suppose $X$ is a finite set and $A\subseteq X$ satisfies $|A|\geq (1-\e)|X|$.  For any partition $\calP$ of $X$, if we let 
$$
\Sigma=\{Y\in \calP: |A\cap Y|\geq (1-a)|Y|\},
$$
 then we have $|\bigcup_{Y\in \Sigma}Y|\geq (1-b)|X|$.
\end{lemma}

\subsection{Weak hypergraph regularity and homogeneity}\label{ss:reg}

In this subsection we introduce weak regularity for $k$-graphs. Since no other form of regularity is used in this paper, we will from here on drop the word ``weak" from all our definitions. Throughout this subsection, $k$ is assumed to be an integer at least $2$. We begin with the definition of $\e$-regular tuples of sets.

\begin{definition}\label{def:regpair}
Suppose $\e>0$, $H=(V,E)$ is a $k$-graph, $V_1,\ldots, V_k\subseteq V$ are nonempty subsets of $V$ (not necessarily distinct).  We say $(V_1,\ldots, V_k)$ is   \emph{$\e$-regular with respect to $H$} if for all $V_1'\subseteq V_1,\ldots, V_k'\subseteq V_k$ satisfying $|V_i'|\geq \e |V_i|$ for each $i\in [k]$, we have
$$
|d_H(V_1,\ldots, V_k)-d_H(V_1',\ldots, V_k')|\leq \e.
$$
\end{definition}

We now define $\e$-regular partitions.

\begin{definition}\label{def:regptn}
Suppose $H=(V,E)$ is a $k$-graph and $\calP$ is a partition of $V$.  We say $\calP$ is \emph{$\e$-regular with respect to $H$} if 
$$
  \sum_{(V_1,\ldots, V_k)\in \Sigma_{\reg}} |V_1\times \cdots \times V_k |\geq (1-\e)|V|^k,
$$
where $\Sigma_{\reg}=\{(V_1,\ldots, V_k)\in \calP^k: (V_1,\ldots, V_k)\text{ is $\e$-regular with respect to $H$}\}$.  In this case, we will also say $\calP$ is \emph{an $\e$-regular partition of $H$}. 
\end{definition}

 The quasirandomness notion appearing in Definition \ref{def:regpair} was first studied independently by Chung \cite{chung} and Haviland-Thomason \cite{HT}.  It was first proved by Chung \cite{chung} that all $k$-graphs admit weak regular partitions  of bounded size \cite{chung}.  The proof is a similar energy increment argument to that typically used to prove Szemer\'{e}di's graph regularity lemma, and produces a similar tower bound (see Theorem \ref{thm:chung} in the introduction).   

We will use the following basic lemmas about regularity. First, large sub-tuples inherit regularity.

\begin{proposition}\label{prop:subsetreg}
Let $0<\e\leq \gamma $, let $H = (V, E)$ be a $k$-graph, and let $(A_1, \ldots, A_k)$ be $\e$-regular with respect to $H$. For each $i\in[k]$ let $A_i'\subseteq A_i$ be such that $|A_i'|\geq \gamma|A_i|.$ Then $(A_1', \ldots, A_k')$ is $2\gamma^{-1}\e$-regular with respect to $H$.
\end{proposition}

We omit the proof of Proposition \ref{prop:subsetreg} as it is essentially identical to the graph  case (see e.g. \cite{KS}). The next lemma says that regular partitions remain somewhat regular under bounded refinements.

\begin{lemma}\label{lem:partrefinement}
 Let $0<\e<1$, let $H = (V, E)$ be a $k$-graph, let $\calP$ be an $\e$-regular partition of $H$, let $\calP'= \{A_1, \ldots, A_\ell\}$ be a partition of $V$ with  $\ell \leq k$, and let $\calQ$ be the common refinement of $\calP$ with $\calP'$. Then $\calQ$ is $(k^2+1)\sqrt{\e}$-regular with respect to $H$.
\end{lemma}
\begin{proof}
    Define
    $$\Sigma_{\reg} = \{(P_1, \ldots, P_k)\in \calP^k: (P_1, \ldots, P_k)\text{ is $\e$-regular with respect to $H$}\}.$$
    Since $\calP$ is $\e$-regular with respect to $H$, we have that 
    $$\left|\bigcup_{(P_1, \ldots, P_k)\in \Sigma_{\reg}}P_{1}\times \ldots\times P_{k}\right|\geq (1-\e)|V|^k.$$
    For all $P\in \calP$ and $1\leq j\leq \ell$, set $P^j=A_j\cap P$. Then $\calQ$ consists of all nonempty sets of the form $P^j$.     We begin by setting notation for the sets $P^j$ from $\calQ$ which are small compared to the original set $P$. Specifically define 
    $$\calQ_{\sml} = \{P^j\in \calQ: |P^j|\leq \sqrt{\e}|P|\}.$$
    We note that 
    \begin{align}\label{al:partrefinement}
    \sum_{P^j\in \calQ_{\sml}}|P^j|\leq \sum_{P^j\in \calQ_{\sml}}\sqrt{\e}|P|\leq \sum_{P\in \calP}k\sqrt{\e}|P| = k\sqrt{\e}|V|.
    \end{align}
    Next we define a set of ``good" $k$-tuples from $\calQ$. In particular, define
    $$\Sigma_\reg' = \{(P_{1}^{j_1}, \ldots, P_{k}^{j_k})\in (\calQ\setminus\calQ_{\sml})^{k}:(P_{1}, \ldots, P_{k})\in \Sigma_{\reg}\}.$$
    By Proposition \ref{prop:subsetreg} every $k$-tuple in $\Sigma_{\reg}'$ is $2\sqrt{\e}$-regular with respect to $H$. Thus it remains to show that $|\bigcup_{(P_{1}^{j_1}, \ldots, P_{k}^{j_k})\in \Sigma_{\reg}'}P_{1}^{j_1}\times \ldots \times P_{k}^{j_k}|\geq (1-(k^2+1)\sqrt{\e})|V|^k$. We have
    \begin{align*}
        \left|\bigcup_{(P_{1}^{j_1}, \ldots, P_{k}^{j_k})\in \Sigma_{\reg}'}P_{1}^{j_1}\times \ldots \times P_{k}^{j_k}\right|&\geq \left|\bigcup_{(P_{1}, \ldots, P_{k})\in \Sigma_{\reg}}P_{1}\times \ldots\times P_{k}\right|-k|V|^{k-1}\sum_{P^j\in \calQ_{\sml}}|P^j|\\
        &\geq (1-\e)|V|^k-k^2\sqrt{\e}|V|^k\\
        &\geq (1-(k^2+1)\sqrt{\e})|V|^k,
    \end{align*}
where the second inequality is by (\ref{al:partrefinement}).
\end{proof}

We now turn to defining homogeneity in $k$-graphs.

\begin{definition}
Suppose $H=(V,E)$ is a $k$-graph and $\e\in (0,1]$.  Given $V_1,\ldots, V_k\subseteq V$, we say $(V_1,\ldots, V_k)$ is \emph{$\e$-homogeneous with respect to $H$} if $d_H(V_1,\ldots, V_k)\in [0,\e)\cup (1-\e,1]$.
\end{definition}

A homogeneous partition is one where most tuples of vertices come from a homogeneous tuple of sets.

\begin{definition}
Suppose $H=(V,E)$ is a $k$-graph and $\e\in (0,1]$. A partition $\calP$ of $V$ is \emph{$\e$-homogeneous with respect to $H$} if 
$$
\sum_{(X_1,\ldots, X_k)\in \Sigma_{\hom}}|X_1\times \cdots \times X_k|\geq (1-\e)|V|^k,
$$
where $\Sigma_{\hom}=\{(X_1,\ldots, X_k)\in \calP^k: d_H(X_1,\ldots, X_k)\in [0,\e)\cup (1-\e,1]\}$.
\end{definition}

Homogeneous tuples are always somewhat regular, as the following exercise shows.

\begin{proposition}\label{prop:homogimpliesreg}
    Let $\e\in (0,1]$ and assume $H=(V,E)$ is a $k$-graph. Suppose $V_1,\ldots, V_k\subseteq V$ satisfy $d_H(V_1, \ldots, V_k)\in [0,\e)\cup (1-\e,1]$. Then $(V_1, \ldots, V_k)$ is $2\e^{1/(k+1)}$-regular with respect to $H$. 
\end{proposition}
\begin{proof}
    Assume $d_H(V_1,\ldots, V_k)\leq \e$ (the other case is similar).  Suppose towards a contradiction $(V_1,\ldots, V_k)$ is not $2\e^{1/(k+1)}$-regular.  Then there are $V_i'\subseteq V_i$ for each $i\in[k]$ such that $|V'_i|\geq 2\e^{1/(k+1)}|V_i|$ and such that 
    $$|d_H(V_1,\ldots, V_k)-d_H(V_1',\ldots, V_k')|> 2\e^{1/(k+1)}.$$
    Since $d_H(V_1, \ldots, V_k)\leq \e$, this implies $d_H(V_1',\ldots, V_k')\geq 2\e^{1/(k+1)}-\e$. But then 
    \begin{align*}
        |\Ebar\cap (V_1\times \ldots\times V_k)|\geq |\Ebar\cap (V_1'\times \ldots\times V_k')|&\geq (2\e^{1/(k+1)}-\e)|V_1'|\ldots|V_k'|\\
        &\geq (2\e^{1/(k+1)}-\e)\e^{k/(k+1)}|V_1|\ldots|V_k|\\
        &=(2\e-\e^{1+k/(k+1)})|V_1|\ldots|V_k|\\
        &> \e|V_1|\ldots|V_k|,
    \end{align*}
    a contradiction.
\end{proof}

\subsection{Hereditary properties}\label{ss:HPs}

In this subsection we introduce definitions and facts related to hereditary $k$-graph properties. 

A \emph{hereditary $k$-graph property} is a non-empty class of finite $k$-graphs closed under isomorphism and induced sub-$k$-graphs. To ease notation at certain points  in the paper, we let $\calG^{(k)}$ denote the hereditary $k$-graph property consisting of \emph{all} finite $k$-graphs.   We now set notation for hereditary $k$-graph properties characterized by  forbidding induced sub-$k$-graphs.  

\begin{definition}
Given $\calF$ a class of finite $k$-graphs, $\Forb(\calF)$ denotes the class of finite $k$-graphs which omit all elements of $\calF$ as induced sub-$k$-graphs.
\end{definition}

It is well known that a class $\calH$ of finite $k$-graphs is a hereditary $k$-graph property if and only if $\calH=\Forb(\calF)$ for some (possibly infinite) $\calF$.  The focus of this paper is on problems which are only interesting when the hereditary property has arbitrarily large elements. For this reason, we set the following convention.

\vspace{2mm}
\noindent{\bf Convention: }   For the rest of the paper, any hereditary $k$-graph property is assumed to contain arbitrarily large $k$-graphs, unless otherwise stated.

\vspace{2mm}

For the sake of completeness, we now repeat the definition of $M_{\calH}$ from the introduction, as this is the main object of study for the paper.

\begin{definition}[Defining $M_{\calH}$]
Suppose $\calH$ is a hereditary $k$-graph property.  Define $M_{\calH}:(0,1]\rightarrow \mathbb{N}$ by setting $M_{\calH}(\e)$ to be the smallest integer $M$ so that all sufficiently large $H$ in $\calH$ admit $\e$-regular partitions with at most $M$ parts.
\end{definition}

Note $M_{\calH}$ is well defined by Theorem \ref{thm:chung}.  Next, we define a well-known notion of ``closeness" for $k$-graphs.

\begin{definition}
    Suppose $H = (V, E)$ and $H' = (V, E')$ are $k$-graphs on the same vertex set. We say $H$ and $H'$ are $\delta$-\emph{close} if $|\Ebar\Delta\Ebar'|\leq \delta|V|^k$.
\end{definition}

Using this, we now define what it means for one hereditary property to be ``close" to another.

\begin{definition}\label{def:closeprop}
Suppose $\calH$ and $\cH'$ are hereditary $k$-graph properties. We say $\calH$ is \emph{close to} $\calH'$ if for all $\delta>0$ there is a natural number $N$ so that for all $H\in \calH$ on at least $N$ vertices, there is $H'\in \cH'$ such that $H$ is $\delta$-close to $H'$. 
\end{definition}
 
Definition \ref{def:closeprop} is relevant in this paper because of the following fact.
 
\begin{fact}\label{fact:close}
Suppose $\calH$ and $\calH'$ are hereditary $k$-graph properties such that $\calH$ is close to $\calH'$. Then for all $0<\e\leq 1/2$, $M_{\calH}(2\e)\leq M_{\calH'}(\e)$.
\end{fact}

We omit the proof of Fact \ref{fact:close} as it is essentially the same as the proof in the graph  case (see Proposition 2.17 in \cite{part2}).  We now define an analogue of $M_{\calH}$ for homogeneous partitions.

\begin{definition}[Defining $M_{\calH}^{hom}$]
Suppose $\calH$ is a hereditary $k$-graph property.  Define $M_{\calH}^{hom}:(0,1]\rightarrow \mathbb{N}\cup\{\infty\}$ by setting $M_{\calH}^{hom}(\e)$ to be the smallest integer $M$ so that all sufficiently large $H$ in $\calH$ admit $\e$-homogeneous partitions with at most $M$ parts, when such an integer exists, and otherwise setting $M_{\calH}^{hom}(\e)=\infty$.
\end{definition}

This function is of interest to us because, in light of Proposition \ref{prop:homogimpliesreg}, we can always roughly upper bound $M_{\calH}$ by $M_{\calH}^{hom}$.

\begin{corollary}
Suppose $\calH$ is a hereditary $k$-graph property.  Then for all $\e\in (0,1]$, $M_{\calH}(\e)\leq M_{\calH}^{hom}((\e/2)^{k+1})$.
\end{corollary}

\subsection{Blowups}\label{ss:blowups}

We now set terminology for turning a $k$-graph  into a larger $k$-graph, by blowing up each vertex into a set of vertices.

\begin{definition}\label{def:blowup}
Suppose $H=(V,E)$ is a $k$-graph and $n\geq 1$ is an integer.  An \emph{$n$-blowup of $H$} is a $k$-graph $H'=(V',E')$ with the following properties.
\begin{enumerate}
\item $V'$ has the form $\bigsqcup_{v\in V}Z_v$, where each $Z_v$ is a set of size $n$;
\item for each $v_1\ldots v_k\in E$, $K_k[Z_{v_1},\ldots, Z_{v_k}]\subseteq E'$ and for each $v_1\ldots v_k\in {V\choose k}\setminus E$, $K_k[Z_{v_1},\ldots, Z_{v_k}]\cap E'=\emptyset$.
\end{enumerate} 
\end{definition}

Note that in the notation of Definition \ref{def:blowup}, an $n$-blowup of $H$ need not be $|V|$-partite. We now define an auxiliary class associated to a hereditary $k$-graph property $\calH$, consisting of all those $k$-graphs with arbitrarily large blowups appearing in $\calH$.

\begin{definition}\label{def:bhgraphs}
 Suppose $k\geq 2$ and $\calH$ is a hereditary $k$-graph property.  Recalling that $\calG^{(k)}$ denotes the class of all finite $k$-graphs, we set
$$
{\bf B}_{\calH}=\{H\in \calG^{(k)}: \text{ for all integers $n\geq 1$, $\calH$ contains an $n$-blowup of $H$}\}.
$$ 
\end{definition}

It was shown in \cite{part2} that ${\bf B}_{\calH}$ is always a hereditary $k$-graph property contained in $\calH$ (we note this implies ${\bf B}_{\calH}$ contains arbitrarily large $k$-graphs, in keeping with our convention).   

\begin{fact}\label{fact:bhhgraphs}
For any hereditary $k$-graph property $\calH$, ${\bf B}_{\calH}$ is a hereditary $k$-graph property satisfying ${\bf B}_{\calH}\subseteq \calH$.   
\end{fact}

One can use ${\bf B}_{\calH}$ to characterize when $\calH$ is close to some $\calH'$, when $\calH'$ is defined via forbidden $k$-graphs.

\begin{fact}\label{fact:blowup}
Suppose $\calH$ is a hereditary $k$-graph property and $\calF$ is a collection of finite $k$-graphs.  The following are equivalent.
\begin{enumerate}
\item $\calH$ is close to $\Forb(\calF)$.
\item ${\bf B}_{\calH}\cap \calF=\emptyset$.  
\end{enumerate}
\end{fact}

We believe Fact \ref{fact:blowup} is well known.  Since we cannot find this exact statement in the literature, we give a brief discussion of how to prove it.  It is an easy exercise to show that if $\calH$ contains an $n$-blowup of $F$ for all $n$, then $\calH$ is not close to $\Forb(\{F\})$. A proof of the $3$-uniform version is included in Theorem A.14 in \cite{part3}, and the $k$-uniform version is essentially identical.  This gives the (1) implies (2) direction. On the other hand, we know from the induced removal lemma for $k$-graphs \cite[Theorem 6]{RS-removal},  that if a $k$-graph $H$ is not close to being $\calF$-free, then it must contain many induced copies of some $F\in \calF$ of bounded size.  By inspecting the proofs of this theorem \cite{RS-removal}, one sees that one can upgrade the conclusion  to obtain a constant sized blowup of $F$ (via the counting lemma).  This strategy can be used to prove the (2) implies (1) direction.   Details of the $3$-uniform case (when $\calF$ is finite) are included in Theorem A.14 of \cite{part3}.   

\section{VC-dimension for hypergraphs and efficient homogeneous partitions}\label{sec:VC}
In this section we define two  versions of VC-dimension for $k$-graphs. It will be useful to first discuss some terminology around set systems.

A \emph{set system} is a pair of the form $(X,\calF)$ where $X$ is the \emph{ground set} and $\calF$ is a family of subsets $\calF\subseteq \calP(X)$.  We will only work with set systems where $X$ is finite. For $Y\subseteq X$, we say $\calF$ \emph{shatters $Y$} if $|\{Y\cap F: F\in \calF\}|=2^{|Y|}$. The VC-dimension of the set system $(X,\calF)$ is then defined to be
$$
\VC(X,\calF)=\max\left(\{0\}\cup \{|Y|: Y\subseteq X\text{ is shattered by $\calF$}\}\right),
$$
 the maximum size of a shattered set.  The VC-dimension of a $k$-graph $G=(V,E)$ is defined as the VC-dimension of a certain neighborhood set system as follows.   

\begin{definition}
For a $k$-graph $G=(V,E)$, the \emph{VC-dimension of $G$}, denoted $\VC(G)$, is the VC-dimension of the set system with ground set $V$ and set collection
$$
\left\{N_G(e): e\in {V\choose k-1}\right\}.
$$
\end{definition}

This notion of VC-dimension has been studied in model theory for many years (see \cite{simon}). Chernikov--Starchenko\cite{CS}, and independently, Fox--Pach--Suk \cite{FPS} showed hypergraphs of bounded VC-dimension admit very small homogeneous partitions. This result will be crucial for us, as Proposition \ref{prop:homogimpliesreg} then implies such $k$-graphs admit very small \emph{regular} partitions as well.  

\begin{theorem} \label{thm:FPS}
For all integers $k\geq 2$ and $d\geq 1$, there are constants $c,c'>0$ such  that the following holds. For all $0<\e<1$, any $k$-graph with VC-dimension at most $d$ has an $\e$-homogeneous partition of size at most $c\e^{-c'}$.
\end{theorem}

We note the constants obtained in \cite{FPS, CS} are explicit, but we omit them here as they are unnecessary for the arguments, and we do not otherwise need to define the dual VC-dimension. It turns out that homogeneous partitions can be guaranteed to exist under weaker assumptions related to  ``slicewise" VC-dimension.  

\begin{definition}[Slice graphs]\label{def:slgraphs}
Given a $k$-uniform hypergraph $H=(V,E)$ and a set $S\in \binom{V}{k-2}$, define the corresponding \emph{slice-graph} to be the graph $H_{S}$ with the following vertex and edge sets. 
$$
V(H_S)=V\text{ and }E(H_S)= \left\{ab\in \binom{V}{2}: S\cup \{a,b\}\in E\right\}.
$$
\end{definition}

We note slice-graphs are usually called link graphs, but we maintain consistency here with \cite{part2} and call them slice graphs.  Definition \ref{def:slgraphs} gives rise to the following version of VC-dimension for $k$-graphs.

\begin{definition}[SVC-dimension]\label{def:svc}
Given a $k$-uniform hypergraph $H=(V,E)$, define its \emph{slicewise VC-dimension} to be
$$
\SVC(H)=\max_{S\in \binom{V}{k-2}}VC(H_S),
$$
where by convention, $\SVC(H)=0$ if $|V(H)|<k-2$.
\end{definition}

It was first shown in \cite{chernikov-towsner} that $k$-graphs of bounded SVC-dimension admit homogeneous partitions of bounded size, although their proof did not provide explicit bounds on the size of such partitions.  A quantitative proof in the $3$-uniform case was first given in \cite{TW}, where a Wowzer type bound was obtained.  Later, it was shown in \cite{part1} that in the $3$-uniform case, the bound can be made  double exponential in $\e^{-1}$. This result was a crucial ingredient in \cite{part2}, which also showed the bound cannot be better than a single exponential.  Subsequent results of Gishboliner, Shapira, and Wigderson show that in all uniformities $k\geq 3$, single exponential bounds are in fact correct.  

\begin{theorem}[Gishboliner, Shapira, and Wigderson \cite{GSW}]\label{thm:gsw}
For all integers $k\geq 3$ and $d\geq 1$, there is a constant $C>0$ so that the following holds.  Suppose $0<\e<1$ and $H$ is a $k$-graph of SVC-dimension at most $d$.  Then there exists an $\e$-homogeneous partition for $H$ with at most $2^{\e^{-C}}$ parts.
\end{theorem}

We now explore what these results tell us about hereditary $k$-graph properties.   We begin by extending the definition of VC-dimension to hereditary $k$-graph properties.  

\begin{definition}[VC-dimension of a hereditary property]
We say a hereditary $k$-graph property $\calH$ has \emph{finite  VC-dimension} if there exists some integer $M$ so that all $H$ in $\calH$ have VC-dimension at most $M$.  In this case we let 
$$
\VC(\calH)=\max\{\VC(H): H\in \calH\}.
$$
Otherwise, we say $\calH$ has \emph{infinite VC-dimension}, and write $\VC(\calH)=\infty$.
\end{definition}

We will need a characterization of $\VC(\calH)<\infty$ in terms of forbidden induced sub-$k$-graphs. We will do this using the following definition.

\begin{definition}\label{def:uk}
    Let $\calU_k(m)$ be the class of all  $k$-graphs $G$ with vertex set of the form 
    $$
   V(G)= \{a^1_i,\ldots, a_i^{k-1}: i\in [m]\}\sqcup \{b_S: S\subseteq [m]\}
    $$
    and edge set $E(G)$ satisfying 
    \begin{align*}
    \text{$i\in S\Rightarrow  a_i^1\ldots a_i^{k-1}b_S\in E(G)$ and $i\in [m]\setminus S\Rightarrow a_i^1\ldots a_i^{k-1}b_S\notin E(G)$.}
    \end{align*}
\end{definition}

We note several facts about $\calU_{k}(m)$. First, we only dictate edges among $k$-tuples of the form $a_i^1\ldots a_i^{k-1}b_S$. We do not say anything about the presence or absence of edges not of the form $a_{i}^1\ldots a_{i}^{k-1}b_S$. Note that Definition \ref{def:uk} implies  $a_{i}^{j_1}\neq a_{i}^{j_2}$ for each $1\leq j_1\neq j_2\leq k-1$ and $1\leq i\leq m$, however it is possible there are other types of overlap among the $a_i^j$ vertices.   

It is an exercise to show that a hereditary $k$-graph property $\calH$ has finite VC-dimension if and only if $\calH$ omits all elements of $\calU_k(m)$ for some integer $m$.

\begin{proposition}\label{prop:vcuk}
Given a hereditary $k$-graph property $\calH$, $\VC(\calH)<\infty$ if and only if there is some $m\geq 1$ such that $\calH\cap \calU_k(m)=\emptyset$.  
\end{proposition}

In light of Fact \ref{fact:close}, we will also be interested in those properties which are merely ``close" to having finite VC-dimension.  

\begin{definition}[Close to finite VC-dimension]
Suppose $\calH$ is a hereditary $k$-graph property. We say $\calH$ is \emph{close to finite VC-dimension} if there is some $\calH'$ with finite VC-dimension so that $\calH$ is close to $\calH'$. Otherwise, we say $\calH$ is \emph{far from finite VC-dimension}.
\end{definition}

Combining with Fact \ref{fact:blowup}, we can characterize the properties close to finite VC-dimension as follows. 

\begin{proposition}\label{prop:vcclose}
Let $\calH$ be a hereditary $k$-graph property. The following are equivalent.
\begin{enumerate}
\item $\calH$ is close to finite VC-dimension.
\item ${\bf B}_{\calH}$ has finite VC-dimension.
\item ${\bf B}_{\calH}\cap \calU_k(m)=\emptyset$ for some $m$.
\end{enumerate}
\end{proposition}
\begin{proof}
The equivalence of (2) and (3) is by Proposition \ref{prop:vcuk}. The equivalence of (1) and (3) is by Proposition \ref{prop:vcuk} and Fact \ref{fact:blowup}. 
\end{proof}

We next define SVC-dimension for hereditary properties. 

\begin{definition}[SVC-dimension of a hereditary property]\label{def:svch}
We say a hereditary $k$-graph property $\calH$ has \emph{finite SVC-dimension} if there exists some integer $M$ so that all $H$ in $\calH$ have SVC-dimension at most $M$.  In this case we let 
$$
\SVC(\calH)=\max\{\SVC(H): H\in \calH\}.
$$
Otherwise, we say $\calH$ has \emph{infinite SVC-dimension}, and write $\SVC(\calH)=\infty$.
\end{definition}

Our next task is to state a characterization of this in terms of forbidden sub-$k$-graphs.  We will use the following notation.

\begin{definition}[Defining $m\otimes_k G$]
Suppose $m\geq 1$ is an integer and $G = ( V, E)$ is a graph. Define $m \otimes_k G$ to be the class of $k$-graphs $H$ with   vertex set $V(H)$ of the form 
$$
V(H)=V\sqcup \left(\bigsqcup_{i=1}^{k-2}\{c_1^i, \ldots, c_m^i\}\right),
$$
where $\{c_j^i: 1\leq i\leq k-2, 1\leq j\leq m\}$ is a set of $m(k-2)$ vertices disjoint from $V$, and with edge set $E(H)$ satisfying 
\begin{align*} 
 \{c^1_{j_1} \ldots c^{k-2}_{j_{k-2}}xy: j_1, \ldots, j_{k-2}\in [m], xy\in E(G)\}\subseteq E(H),
 \end{align*}
 and 
 \begin{align*} 
 \left\{c^1_{j_1} \ldots c^{k-2}_{j_{k-2}}xy: j_1, \ldots, j_{k-2}\in [m], xy\in {V\choose 2}\setminus E(G)\right\}\cap E(H)=\emptyset.
 \end{align*}
\end{definition}

We extend the above notation to classes of graphs $G$ as follows.

\begin{definition}
Given integers $k\geq 2$ and $m\geq 1$ and a class $\calC$ of finite graphs, set $m\otimes_k\calC=\bigcup_{G\in \calC}m\otimes_kG$.
\end{definition}

This allows us to easily characterize hereditary $k$-graph properties with finite SVC-dimension.

\begin{proposition}\label{prop:svcuk}
Suppose $\calH$ is a hereditary $k$-graph property.  Then $\SVC(\calH)<\infty$ if and only if there is $m\geq 1$ such that $\calH\cap (1\otimes_k \calU_2(m))=\emptyset$.  
\end{proposition}

As above, it will turn out that the more fundamental notion for this paper is that of being ``close" to a property with finite SVC-dimension.

\begin{definition}[Close to finite SVC-dimension]
Suppose $\calH$ is a hereditary $k$-graph property. We say $\calH$ is \emph{close to finite SVC-dimension} if there is some $\calH'$ with finite SVC-dimension so that $\calH$ is close to $\calH'$.  Otherwise, we say $\calH$ is \emph{far from finite SVC-dimension}.
\end{definition}

Combining with Fact \ref{fact:blowup}, we can characterize the properties close to finite SVC-dimension as follows. 

\begin{proposition}\label{prop:svcequiv}
Let $\calH$ be a hereditary $k$-graph property. The following are equivalent.
\begin{enumerate}
\item $\calH$ is close to finite SVC-dimension.
\item ${\bf B}_{\calH}$ has finite SVC-dimension.
\item ${\bf B}_{\calH}\cap(1\otimes_k \calU_2(m))=\emptyset$ for some $m$.
\item $\calH\cap (m\otimes_k\calU_2(m))=\emptyset$ for some $m$.
\end{enumerate}
On the other hand, if (1)-(4) fail, then one of the following holds.
\begin{enumerate}[(a)]
\item $(n\otimes_k G)\cap \calH\neq \emptyset$ for every $n\geq 1$ and finite bipartite graph $G$.
\item $(n\otimes_k G)\cap \calH\neq \emptyset$ for every $n\geq 1$ and finite cobipartite graph $G$.
\item $(n\otimes_k G)\cap \calH\neq \emptyset$ for every $n\geq 1$ and finite split graph $G$.
\end{enumerate}
\end{proposition}
\begin{proof}
The equivalence of (2) and (3) is by Proposition \ref{prop:svcuk}. The equivalence of (1) and (3) is by Proposition \ref{prop:svcuk} and Fact \ref{fact:blowup}.  The equivalence of (3) and (4) is an exercise.

Assume now (1)-(4) fail.  The fact that one of (a)-(c) holds can be proved via a similar strategy as in the 3-uniform case (see Theorem 4.18 of \cite{part2}), which uses the fact that a hereditary \emph{graph} property has infinite VC-dimension if and only if it contains every bipartite, every split, or every cobipartite graph.  
\end{proof}

It was shown in \cite{TW} that in the 3-uniform case, a hereditary property admits homogeneous partitions if and only if it is close to finite SVC-dimension.  Using Theorem \ref{thm:gsw} and similar arguments to \cite{TW} (see the proof of Theorem 2.30 there), one can extend this characterization to arity $k$.   

\begin{theorem} 
Suppose $\calH$ is a hereditary $k$-graph property. The following are equivalent.
\begin{enumerate}
\item $\calH$ is close to finite SVC-dimension.
\item For all $\e\in (0,1)$, there exists an integer $M=M(\e)$ so that every $k$-graph in $\calH$ admits an $\e$-homogeneous partition of size at most $M$.
\end{enumerate} 
\end{theorem}

In other words, $M_{\calH}^{\hom}<\infty$ if and only if $\calH$ is close to finite $\SVC$-dimension.

\section{Jump from Exponential to Tower}
In this section we prove that there is a jump between exponential growth and tower growth. We begin by stating a result by Fox and Lov\'{a}sz from \cite{FL}, which tells us there exist bipartite graphs for which any regular partition requires tower-many parts.  As in \cite{part2}, we use the fact that an inspection of \cite{FL} shows it can be easily adapted to build instead a cobipartite or split example.

\begin{theorem}[Fox, Lov\'{a}sz \cite{FL}]\label{thm:FL}
There is a constant $c>0$ so that for all sufficiently small $\e>0$ and all $n_0$, there exist $n\geq n_0$, graphs $G_1,G_2,G_3$ on the same vertex set of size $2n$,  such that $G_1$ is bipartite, $G_2$ is cobipartite, and $G_3$ is a split graph, and such that every $\e$-regular partition of $G_1$, $G_2$, or $G_3$ has size at least $Tw(c\e^{-2})$.
\end{theorem}

We now prove an analogue of Lemma 5.2  from \cite{part2}. It says that if $H \in n\otimes_k G$ and $H$ has a regular partition, then $G$ has a regular partition of a similar size.  

\begin{lemma}\label{lem:ntimesG}
 There is a constant $D_k$ such that for all $\e>0$ sufficiently small, the following holds.   Suppose $G$ is a graph with $2n$ vertices and $H \in n \otimes_k G$. If $H$ has an $\epsilon$-regular partition of size $t$ then $G$ has a $D_k\e^{1/2}$-regular partition of size at most $2t +2$.
\end{lemma}
\begin{proof}
We proceed via an argument very similar to Lemma 5.2 from \cite{part2}.  By assumption, $H$ has vertex set of the form 
$$
C_1\sqcup \ldots \sqcup C_{k-2}\sqcup V(G),
$$
where $|C_1|=\ldots=|C_{k-2}|=n$, and edge set $E(H)$ satisfying
\begin{align}\label{al:eg}
& \{c_1\ldots c_{k-2}xy: c_1\in C_1,\ldots, c_{k-2}\in C_{k-2}, xy\in E(G)\}\subseteq E(H),\text{and}\\
\nonumber  &\left\{c_1\ldots c_{k-2}xy: c_1\in C_1,\ldots, c_{k-2}\in C_{k-2}, xy\in {V(G)\choose 2}\setminus E(G)\right\}\cap E(H)=\emptyset.
 \end{align}
Let $\calP$ be an $\e$-regular partition of $H$ with size $t$, and let $\calQ$ be the common refinement of $\calP$ with $\{C_1,\ldots, C_{k-2},V(G)\}$.  By Lemma \ref{lem:partrefinement}, $\calQ$ is $d_k\e^{1/2}$-regular where $d_k=k^2+1$. Let 
$$
\Sigma_{\reg}=\{(U_1,\ldots, U_k)\in \calQ^k: (U_1,\ldots, U_k)\text{ is }d_k\e^{1/2}\text{-regular with respect to $H$}\},
$$ 
and let 
$$
\calE_{\reg}=\bigcup_{(U_1,\ldots, U_k)\in \Sigma_{\reg}}U_1\times \cdots \times U_k.
$$
By construction, $|\calE_{\reg}|\geq (1-d_k\e^{1/2})|V(H)|^k$.  Since $|V(G)|=2n$ and each $C_i$ has size $n$, we have
\begin{align*}
|\calE_{\reg}\cap (V(G)\times V(G)\times C_1\times \cdots \times C_{k-2})|&\geq |V(G)\times V(G)\times C_1\times \cdots \times C_{k-2}|-d_k\e^{1/2}(kn)^k\\
&\geq (1-k^kd_k\e^{1/2})|V(G)\times V(G)\times C_1\times \cdots \times C_{k-2}|.
\end{align*}
Since $V(G)\times V(G)\times C_1\times \cdots \times C_{k-2}$ is a union of sets of the form $X_1\times \cdots \times X_k$ for $X_1,\ldots, X_k\in \calQ$, simple averaging implies there exist some $Z_1,\ldots Z_{k-2}\in \calQ$ with $Z_i\subseteq C_i$ for each $i\in [k-2]$, such that if we set 
$$
\Theta=\{(X,Y)\in \calQ^2: X,Y\subseteq V(G)\text{ and }X\times Y\times Z_1\times \cdots \times Z_{k-2}\subseteq \calE_{\reg}\},
$$
then 
$$
\left|\bigcup_{(X,Y)\in \Theta}X\times Y\right|\geq (1-k^kd_k\e^{1/2})|V(G)|^2.
$$
For all nonempty $A,B\subseteq V(G)$, (\ref{al:eg}) implies 
$$
|\overline{E(H)}\cap (A\times B\times Z_1\times \cdots\times Z_{k-2})|=|\overline{E(G)}\cap (A\times B)||Z_1|\cdots |Z_{k-2}|,
$$
and thus $d_G(A,B)=d_H(A,B,Z_1,\ldots, Z_{k-2})$. Consequently, if $(A,B)\in \Theta$, then for all $A'\subseteq A$ and $B'\subseteq B$ satisfying $|A'|\geq d_k\e^{1/2}|A|$ and $|B'|\geq d_k\e^{1/2} |B|$, we have 
$$
\left|d_G(A,B)-d_G(A',B')\right|=\left|d_H(A,B,Z_1,\ldots, Z_{k-2})-d_H(A',B',Z_1,\ldots, Z_{k-2})\right|\leq d_k\e^{1/2}. 
$$
We have now shown every $(X,Y)\in \Theta$ is $k^kd_k\e^{1/2}$-regular in $G$, which concludes the proof.  
 \end{proof}

We now show that if $\calH$ is far from finite $\SVC$-dimension, then $M_{\calH}$ is bounded below by a tower type function.

\begin{theorem}\label{thm:twr}
Let $\calH$ be a hereditary $k$-graph property. If $\calH$ is far from finite SVC-dimension, then for all sufficiently small $\e>0$, $M_{\calH}(\e)\geq Tw(\e^{-C})$ for some $C>0$. 
\end{theorem}
\begin{proof}
Let $\calH$ be a hereditary $k$-graph property and suppose $\cH$ is far from finite SVC-dimension.  By Proposition \ref{prop:svcequiv}, one of the following holds. 
\begin{enumerate}[(a)]
\item $(n\otimes_k G)\cap \calH\neq \emptyset$ for every $n\geq 1$ and finite bipartite graph $G$.
\item $(n\otimes_k G)\cap \calH\neq \emptyset$ for every $n\geq 1$ and finite cobipartite graph $G$.
\item $(n\otimes_k G)\cap \calH\neq \emptyset$ for every $n\geq 1$ and finite split graph $G$.
\end{enumerate}
 Suppose $\e>0$ is sufficiently small and let $n_0$ be sufficiently large. Let $D_k$ be as in Lemma \ref{lem:ntimesG}, and let $G$ and $n\geq n_0$ be as in Theorem \ref{thm:FL} with parameter $D_k\e^{1/2}$, with $G$ chosen to be bipartite, cobipartite, or split, depending on whether (a), (b), or (c) hold, respectively. Then there is some $H\in \cH\cap (n \otimes_k G)$. Let $t_H$ be the size of the smallest $\e$-regular partition of $H$, and let $t_G$ be the size of the smallest $D_k\e^{1/2}$-regular partition of $G$. Then by Lemma \ref{lem:ntimesG}, $t_G\leq 2t_H+2$, and by Theorem \ref{thm:FL}, $t_G\geq Tw(cD_k^{-2}\e^{-1})$ for some constant $c$. Thus there is a constant $C'$ such that $M_{\cH}(\e)\geq t_H\geq \frac{t_G-2}{2}\geq Tw(\e^{-C'})$. 
\end{proof}

We are now ready to characterize which hereditary $k$-graph properties require tower size regular partitions.

\begin{corollary}\label{cor:tower}
Suppose $k\geq 3$ is an integer and $\calH$ is a hereditary $k$-graph property. Then one of the following holds.
\begin{enumerate}
\item $\calH$ is close to finite SVC-dimension. In this case there is $C>0$ such that for all sufficiently small $\e>0$, $M_{\calH}(\e)\leq 2^{\e^{-C}}$.
\item $\calH$ is far from finite SVC-dimension. In this case there is $C>0$ such that for all sufficiently small $\e>0$, $M_{\calH}(\e)\geq Tw(\e^{-C})$. 
\end{enumerate}
\end{corollary}
\begin{proof}
    First, suppose $\calH$ is close to finite SVC-dimension. Then there is $\calH'$ with finite SVC-dimension such that $\cH$ is close to $\cH'$. By Theorem \ref{thm:gsw} and Proposition \ref{prop:homogimpliesreg}, there is a constant $C'$ such that $M_{\cH'}(\e)\leq 2^{\e^{-C'}}$. Therefore, by Fact \ref{fact:close}, there is a constant $C$ such that $M_{\cH}(\e)\leq 2^{\e^{-C}}$. On the other hand, if $\cH$ is far from finite $\SVC$-dimension, then by Theorem \ref{thm:twr}, there is a constant $C>0$ such that $M_\cH(\e)\geq Tw(\e^{-C})$.
\end{proof}

\section{A lower bound lemma}\label{sec:LB}

This section proves a sufficient condition for producing lower bounds from blowups.  This result, Lemma \ref{lem:blowup}, will be the key ingredient in our constant, polynomial, and exponential lower bounds.   Stating and proving this lemma will require a few preliminaries, in close analogy to \cite[Section 6]{part2}.  We begin with a definition of when two vertices are ``twins" in a $k$-graph. This notion has appeared many times in the literature (e.g. \cite[Section~2]{HM88}). 

\begin{definition}[Twins]\label{def:twins}
Suppose $H=(V,E)$ is a $k$-graph and $x,x'$ are elements of $V$.  We say $x,x'$ are \emph{twins in $H$}, denoted  $x\sim_Hx'$, if and only if
$$
N_H(x)\cap {V\setminus \{x,x'\}\choose k-1}=N_H(x')\cap {V\setminus \{x,x'\}\choose k-1}.
$$
\end{definition}

It is not too difficult to show that $\sim_H$ is always an equivalence relation on the vertices of $H$, the classes of which are cliques or independent sets.  We refer to these classes as the $\sim_H$-classes in $H$.  In this section, we will work with the following variant of Definition \ref{def:twins} relativized to subsets. This is the natural $k$-ary analogue of Definition 6.2 in \cite{part2}.

\begin{definition}[Local twins]\label{def:localtwin}
Suppose $H=(V,E)$ is a $k$-graph, $X,Y\subseteq V$ are nonempty subsets satisfying either $X\cap Y=\emptyset$  or $X=Y$, and $x,x'$ are elements of $X$.  We say $x$ and $x'$ are \emph{$Y$-twins in $H$}, denoted $x\sim_H^Y x'$, if and only if 
$$
N_H(x)\cap {Y\setminus \{x,x'\}\choose k-1}=N_H(x')\cap {Y\setminus \{x,x'\}\choose k-1}.
$$
\end{definition}

It is an exercise  to check $\sim_H^Y$ is an equivalence relation on the set $X$ (this uses the assumption that $X$ and $Y$ are either disjoint or equal).  We refer to its equivalence classes as the $\sim_H^Y$-classes in $X$.  We now recall  some notation set out in Subsection \ref{ss:not}. Given a $k$-graph $H=(V,E)$, we let $E^1=E$, $E^0={V\choose k}\setminus E$, $H^1=H$, and $H^0=(V,E^0)$. Given any $S\in {V\choose k-1}$, this always yields a partition of $V$ of the form
$$
V=N_{H^1}(S)\cup N_{H^0}(S)\cup S.
$$
We next state the $k$-graph analogue of Lemma 6.4 in \cite{part2}, which tells us certain combinations of sets like those in the partition above are constrained to single equivalence classes, in the sense of Definition \ref{def:localtwin}.  For the sake of completeness, we include the proof below, but we note that it is basically identical to the proof of Lemma 6.4 in \cite{part2}.  

\begin{lemma}\label{lem:simclasslemma}
Suppose $H=(V,E)$ is a $k$-graph and $X,Y$ are nonempty subsets of $V$ which are either disjoint or equal to each other.  For any function $\alpha:{Y\choose k-1}\rightarrow \{0,1\}$, the following set is contained in a single $\sim_H^Y$ class in $X$:
$$
X\cap \bigcap_{S\in {Y\choose k-1}}\left(S\cup N_{H^{\alpha(S)}}(S)\right).
$$
\end{lemma}
\begin{proof}
Fix $\alpha:{Y\choose k-1}\rightarrow \{0,1\}$, and to ease notation, set 
$$
W=X\cap  \bigcap_{S\in {Y\choose k-1}}(S\cup N_{H^{\alpha(S)}}(S)).
$$
  Suppose towards a contradiction $W$ is not a subset of a single $\sim^Y_H$-class in $X$.  Then there exist $x,x'\in W$  such that $x\nsim^Y_H x'$.  By definition of $\sim^Y_H$, there is $S\in {Y\setminus\{x,x'\}\choose k-1}$ and $\tau\in \{0,1\}$   such that $S\cup \{x\}\in E^{\tau}$ and $S\cup \{x'\}\in E^{1-\tau}$. Since $x,x'\in W$, we have $\{x,x'\}\subseteq S\cup N_{H^{\alpha(S)}}(S)$. Since $S\cap \{x,x'\}=\emptyset$, this implies $\{x,x'\}\subseteq   N_{H^{\alpha(S)}}(S)$. However, this is impossible, since either  $\alpha(S)=\tau$ and $x'\notin N_{H^{\alpha(S)}}(S)$, or $\alpha(S)=1-\tau$ and $x\notin N_{H^{\alpha(S)}}(S)$.
\end{proof}

We next define a slightly more complicated notion of blowup, where we keep track of the sizes of some of the blown up vertex sets.

\begin{definition}
    Let $H = (V, E_H)$ be a $k$-graph, let $X,Y\subseteq V$ be sets, and let $n_1,n_2\geq 1$ be integers. An \emph{$(X,Y; n_1,n_2)$-blowup of $H$} is any $k$-graph $\mathbf{H}$ with the following properties. 
    \begin{itemize}
        \item $\mathbf{H}$ has vertex set of the form
        $$
        V(\mathbf{H}) = \bigsqcup_{a\in V}Z_a,
        $$
        where for each $a\in X$, $|Z_a|=n_1$ and for each $b\in Y$, $|Z_b|=n_2$. 
        \item $\mathbf{H}$ has edge set satisfying 
        $$\bigcup_{a_1\ldots a_{k}\in E(H)}K_k[Z_{a_1}, \ldots , Z_{a_{k}}]\subseteq E(\mathbf{H})$$
        and 
        $$\left(\bigcup_{a_1\ldots a_{k}\in {V\choose k}\setminus E(H)}K_k[Z_{a_1}, \ldots , Z_{a_{k}}]\right)\cap E(\mathbf{H}) = \emptyset$$
    \end{itemize}
\end{definition}

We are now ready to state and prove the main result of this section.  The reader may wish to review Lemma \ref{lem:averaging} as well as the notation directly preceding it.

\begin{lemma}\label{lem:blowup}
Let $H = (V, E_H)$ be a $k$-graph, let $n_1,n_2\geq 1$ be integers, and let $X,Y\subseteq V$ be subsets satisfying $|X|n_1=|Y|n_2$ and either $X\sqcup Y=V$ or $X=Y=V$.  Assume   
    $$
    0<\e<\left(\frac{1}{2(k^2+1)|Y|}\right)^{2^{k}}.
    $$
   Suppose $\mathbf{H}$ is an $(X,Y;n_1,n_2)$-blowup of $H$ with vertex set 
    $V(\mathbf{H}) = \bigcup_{u\in V}Z_u$.  For any $\e$-regular partition $\calP$ of $\mathbf{H}$, there is a subset $\calP'\subseteq \calP$ such that the following hold.
    \begin{enumerate}[(i)]
        \item  $|\bigcup_{P\in \calP'}P\cap {\bf X}|\geq (1-2(k^2+1)\e^{1/2^k})|{\bf X}|$, where ${\bf X} = \bigcup_{u\in X}Z_u$.
        \item For all $P\in \calP'$, there is a single $\sim_H^Y$-class $X'$ in $X$  such that 
        $$
        |P\cap {\bf X}'|\geq (1- \e^{1/4})|P\cap {\bf X}|,
        $$
         where ${\bf X}' = \bigcup_{u\in X'}Z_u$ and ${\bf X} = \bigcup_{u\in X}Z_u$.
    \end{enumerate}
\end{lemma}
\begin{proof}
To ease notation throughout the proof, we let ${\bf X} = \bigcup_{u\in X}Z_u$ and ${\bf Y}=\bigcup_{u\in Y}Z_u$. Our assumptions imply $n_1|X|=n_2|Y|$, and either ${\bf X}={\bf Y}=V({\bf H})$  or ${\bf X}\sqcup {\bf Y}=V({\bf H})$.  Combined with the definition of ${\bf H}$ being an $(X,Y; n_1,n_2)$-blowup, this yields $|{\bf X}|=|{\bf Y}|$ and $|V({\bf H})|\leq 2|{\bf X}|=2|{\bf Y}|$.

   Let $\calQ$ be  the common refinement of $\calP$ and the partition of $V({\bf H})$ generated by the sets ${\bf X},{\bf Y}$.  Note that because either ${\bf X}={\bf Y}=V({\bf H})$ or ${\bf X}\sqcup {\bf Y}=V({\bf H})$, every set in $\calQ$ has the form $P\cap {\bf X}$ or $P\cap {\bf Y}$ for some $P\in \calP$.
   
By Lemma \ref{lem:partrefinement}, $\calQ$ is $\eta$-regular with respect to ${\bf H}$, where $\eta:=(k^2+1)\sqrt{\e}$. We now set notation for the set of $\eta$-regular $k$-tuples from $\calQ$ and the corresponding set of underlying vertex tuples. In particular, set
   \begin{align*}
   \Sigma_{\reg}  &= \{(Q_1, \ldots, Q_k)\in \calQ^k: (Q_1, \ldots, Q_k)\text{ is $\eta$-regular with respect to $\mathbf{H}$}\}\text{ and }\\
   \calE_{\reg}&=\bigcup_{(Q_1,\ldots, Q_k)\in \Sigma_{\reg}}Q_1\times \cdots \times Q_k.
   \end{align*}
 Since  $\calQ$ is $\eta$-regular with respect to ${\bf H}$, we know $|\calE_{\reg}|\geq (1-\eta)|V(\mathbf{H})|^k$.
     
 Next, we inductively define,  for each $0\leq i\leq k-1$, a ``good" collection $\Sigma_{\reg}^{(k-i)}$ of  tuples $(Q_1,\ldots, Q_{k-i})$ from $\calQ^{k-i}$, along with the corresponding set of underlying vertex tuples $\calE^{(k-i)}_{\reg} \subseteq V({\bf H})^{k-i}$. As we proceed, we will ensure at each step that  $\calE_{\reg}^{(k-i)}$ is large, in particular, that $|\calE_{\reg}^{(k-i)}|\geq (1-\eta^{1/2^{i}})|V({\bf H})|^{k-i}$. 
    
    \begin{itemize}
    \item Base case: Set $\Sigma_{\reg}^{(k)}=\Sigma_{\reg}$ and $\calE^{(k)}_{\reg} = \calE_{\reg}$.  By above, we know 
    $$
    |\calE_{\reg}^{(k)}|\geq (1-\eta)|V(\mathbf{H})|^k.
    $$
    \item Step $i$:   Suppose  $1\leq i\leq k-1$, and assume by induction we have defined a set $\Sigma_{\reg}^{(k-i+1)}\subseteq \calQ^{k-i+1}$ and $\calE_{\reg}^{(k-i+1)}=\bigcup_{(Q_1,\ldots, Q_{k-i+1})\in \Sigma_{\reg}^{(k-i+1)}}Q_1\times \cdots \times Q_{k-i+1}$ satisfying $|\calE_{\reg}^{(k-i+1)}|\geq (1-\eta^{1/2^{i-1}})|V(\mathbf{H})|^{k-i+1}$. We then let  
    $$
    \calE_{\reg}^{(k-i)} = \{\vbar\in V(\mathbf{H})^{k-i}: |N_{\calE_{\reg}^{(k-i+1)}}(\vbar)|\geq (1-\eta^{1/2^{i}})|V(\mathbf{H})|\}.
    $$
    By Lemma \ref{lem:averaging} $|\calE_{\reg}^{(k-i)}|\geq (1-\eta^{1/2^{i}})|V(\mathbf{H})|^{k-i}$.  Since $\calE_{\reg}^{(k-i+1)}$ is a disjoint union of sets of the form $Q_1\times \cdots \times Q_{k-i+1}$ with $Q_1,\ldots, Q_{k-i+1}\in \calQ$, there is a set $\Sigma_{\reg}^{(k-i)}\subseteq \calQ^{k-i}$ such that $\calE_{\reg}^{(k-i)}=\bigcup_{(Q_1,\ldots, Q_{k-i})\in \Sigma_{\reg}^{(k-i)}}Q_1\times \cdots\times Q_{k-i}$. This completes step $i$.
    \end{itemize}
\vspace{2mm}

At the end of this construction we obtain $\calE_{\reg}^{(1)}\subseteq V({\bf H})$ satisfying  $|\calE_{\reg}^{(1)}|\geq (1-\eta^{1/2^{ k-1}})|V(\mathbf{H})|$. Moreover, by construction, $\calE_{\reg}^{(1)}$ is a union of sets from $\calQ$.  Since $\calQ$ is the common refinement of $\calP$ with the partition generated by $\{{\bf X}, {\bf Y}\}$,  there is a set $\calP'\subseteq \calP$ such that
$$
\bigcup_{P\in \calP'}P\cap {\bf X}= \calE_{\reg}^{(1)}\cap {\bf X}.
$$
We will prove $\calP'$ has the desired properties appearing in conditions (i) and (ii) of Lemma \ref{lem:blowup}.  First, by our lower bound on $|\calE_{\reg}^{(1)}|$ and the fact that $|{\bf X}| \geq \frac{1}{2}|V(\mathbf{H})|$ we have that  
    $$
   \left|\bigcup_{P\in \calP'}P\cap {\bf X}\right|=  |\calE_{\reg}^{(1)}\cap {\bf X}|\geq |{\bf X}|-\eta^{1/2^{k-1}}|V(\mathbf{H})|\geq |{\bf X}|-2\eta^{1/2^{k-1}}|{\bf X}|\geq (1-2(k^2+1)\e^{1/2^{k}})|{\bf X}|,
    $$
where the last inequality uses the definition of $\eta$. We have now shown (i) holds.  

The rest of the proof is devoted to the proof of (ii).   To this end, fix $P\in \calP'$ for the rest of the proof. Our goal is to show  there is a single $\sim_H^Y$-class $X'$ in $X$  such that, setting ${\bf X}' = \bigcup_{u\in X'}Z_u$, we have $|P\cap {\bf X}'|\geq (1- \e^{1/4})|P\cap {\bf X}|$.
         
Recall that any $S=\{a_2,\ldots, a_{k}\}\in {Y\choose k-1}$ induces a natural partition of ${\bf X} \setminus (Z_{a_2}\cup \ldots\cup Z_{a_k})$ given by 
    $$
    {\bf X}\setminus (Z_{a_2}\cup \ldots\cup Z_{a_k})= \left(\bigcup_{u\in X\cap N_{H^1}(S)}Z_u\right)\sqcup \left(\bigcup_{u\in X\cap N_{H^0}(S)}Z_u\right).
    $$
 These partitions will play an important role in the rest of the proof. To ease notation, for each $S\in {Y\choose k-1}$ and  $\alpha\in \{0, 1\}$ we let
    $$
    B_S^\alpha = \bigcup_{u\in X\cap N_{H^\alpha}(S)}Z_u.
    $$
    Note any such $B_S^{\alpha}$ is by definition a subset of ${\bf X}$.  We next show that $P$ cannot substantially intersect both $B_S^0$ and $B_S^1$ for any $S\in {Y\choose k-1}$.  More formally we prove the following claim.
    
    \begin{claim}\label{cl:size_claim}
    For all $S\in {Y\choose k-1}$, $\min_{\alpha\in \{0,1\}}\left|P\cap B_S^\alpha\right|\leq \eta |P\cap {\bf X}|$.
    \end{claim}
    \begin{proof}[Proof of Claim]
        Fix  $S=\{a_2,\ldots, a_k\}\in {Y\choose k-1}$.  By construction, there exists some $Q_1\in \calQ$ such that $P\cap {\bf X}=Q_1$.   By definition of $\calP'$, we have $Q_1\subseteq \calE_{\reg}^{(1)}\cap {\bf X}$. If $Q_1\cap B_S^0=\emptyset$ or $Q_1\cap B_S^1=\emptyset$, then we are done.  So assume $Q_1\cap B_S^0$ and $Q_1\cap B_S^1$ are both nonempty.
        
We will inductively construct a tuple $(Q_2, \ldots, Q_k)$ such that $(Q_1, Q_2,\ldots, Q_k)\in \Sigma_{\reg}^{(k)}$,  and such that for each $i\in\{2, \ldots, k\}$,
        $$
        |Q_i\cap Z_{a_i}|\geq \eta |Q_i|.
        $$
        Suppose $2\leq i\leq k$, and assume by induction we have defined $Q_j$ for each $1\leq j\leq i-1$, such that $(Q_1,\ldots, Q_{i-1})\in \Sigma_{\reg}^{(i-1)}$, and such that for all $2\leq j\leq i-1$, $|Q_j\cap Z_{a_j}|\geq \eta |Q_j|$ (note this is vacuously true in our base case, $i=2$).   Let 
        $$
        \calQ(i)=\{W\in \calQ: W\subseteq {\bf Y}\text{ and }(Q_1, \ldots, Q_{i-1}, W)\in \Sigma_{\reg}^{(i)}\}.
        $$
         Since $(Q_1,\ldots, Q_{i-1})\in \Sigma_{\reg}^{(i-1)}$ and $Z_{a_i}\subseteq {\bf Y}$, we have that the following holds, where $(q_1,\ldots, q_{i-1})$ is any element of $Q_1\times \cdots \times Q_{i-1}$: 
        \begin{align*}
            \left|Z_{a_i}\cap {\bf Y}\cap \left(\bigcup_{W\in \calQ(i)}W\right)\right| =\left|Z_{a_i}\cap  \left(\bigcup_{W\in \calQ(i)}W\right)\right|&=  \left|Z_{a_i}\cap \left(N_{\calE_{\reg}^{(i)}}(q_1,\ldots, q_{i-1}) \right)\right| \\
            &\geq \left|Z_{a_i}\right|- \eta^{1/2^{k-i+1}}|V(\mathbf{H})|\\
            &= | {\bf Y}||Y|^{-1}- \eta^{1/2^{k-i+1}}|V({\bf H})|\\ 
            &\geq ((2|Y|)^{-1}- \eta^{1/2^{k-i+1}})|V({\bf H})|\\
            &\geq \eta |V({\bf H})|,
        \end{align*}
        where the last inequality uses the definition of $\eta$ and the assumed upper bound on $\e$.  Thus
        \begin{align*}
            \left|Z_{a_i}\cap {\bf Y}\cap \left(\bigcup_{W\in \calQ(i)}W\right)\right|&\geq \eta |{\bf Y}|\geq \eta \left|{\bf Y}\cap \left(\bigcup_{W\in \calQ(i)}W\right)\right|.
        \end{align*}
 Therefore by pigeonhole principle there is $Q_i\in \calQ(i)$ such that $|Z_{a_i}\cap Q_i|\geq \eta |Q_i|$. By definition of $\calQ(i)$, we have $(Q_1, \ldots, Q_i)\in \Sigma_{\reg}^{(i)}$.

 After $k$ steps, this constructs a sequence $(Q_1, Q_2,\ldots, Q_k)\in \Sigma_{\reg}^{(k)}$, where $Q_2,\ldots, Q_k\subseteq {\bf Y}$, and such that for each $i\in \{2, \ldots, k\}$ we have $|Q_i\cap Z_{a_i}|\geq \eta |Q_i|$. Note that
    $$
    d_{\mathbf{H}}(Q_1\cap B_S^1, Q_2\cap Z_{a_2}, \ldots, Q_k\cap Z_{a_k}) = 1
    $$
    and
    $$
    d_{\mathbf{H}}(Q_1\cap B_S^0, Q_2\cap Z_{a_2}, \ldots, Q_k\cap Z_{a_k}) = 0.
    $$
    Therefore, our upper bound on $\e$ and the definition of $\eta$ imply that either
    $$
    |d_{\mathbf{H}}(Q_1\cap B_S^1, Q_2\cap Z_{a_2}, \ldots, Q_k\cap Z_{a_k})-d_{\mathbf{H}}(Q_1, \ldots, Q_k)|>\eta
    $$
    or
    $$
    |d_{\mathbf{H}}(Q_1\cap B_S^0, Q_2\cap Z_{a_2}, \ldots, Q_k\cap Z_{a_k})-d_{\mathbf{H}}(Q_1, \ldots, Q_k)|>\eta.
    $$
    By construction, $(Q_1,\ldots, Q_k)\in \Sigma_{\reg}^{(k)}=\Sigma_{\reg}$ and $|Q_i\cap Z_{a_i}| \geq \eta |Q_i|$ for all $i\in \{2, \ldots, k\}$. Consequently, we must have that either
    $$
    |Q_1\cap B_S^1|\leq \eta |Q_1|\text{ or }|Q_1\cap B_S^0|\leq \eta |Q_1|.
    $$
    This completes our proof of the claim.
    \end{proof}
   
   We now finish our proof that (ii) holds.  For each $S\in {Y\choose k-1}$, let $\alpha(S)$ be such that $|P\cap B_S^{1-\alpha(S)}|\leq \eta |P\cap {\bf X}|$ (this exists by Claim \ref{cl:size_claim}).   By Lemma \ref{lem:simclasslemma}, there exists a single $\sim_H^Y$-class $X'$ containing $X\cap \bigcap_{S\in {Y\choose k-1}}(S\cup N_{H^{\alpha(S)}}(S))$. Consequently, if we set ${\bf X}'=\bigcup_{u\in X'}Z_u$, we have 
   $$
{\bf X}\cap \bigcap_{S\in {Y\choose k-1}}\left(\bigcup_{u\in (S\cup N_{H^{\a(S)}}(S))}Z_u\right)\subseteq {\bf X}'.
 $$ 
Thus
\begin{align*}
 |P\cap {\bf X}'|\geq   \left|P\cap {\bf X}\cap \bigcap_{S\in {Y\choose k-1}}\left(\bigcup_{u\in (S\cup N_{H^{\a(S)}}(S))}Z_u\right)\right|&\geq |P\cap {\bf X}|-|Y|^{k-1}\eta|P\cap {\bf X}|\\
 &\geq (1- \e^{1/4})|P\cap {\bf X}|,
 \end{align*}
    where the last inequality is by definition of $\eta$ and our upper bound on the size of $\e$.   This concludes our proof of (ii).
   
\end{proof}

We now deduce a corollary which will be convenient in our applications, in close analogy to Corollary 6.6 in \cite{part2}.

\begin{corollary}\label{cor:blowup}
 Let $H $ be a $k$-graph, let $n_1,n_2\geq 1$ be integers, and let $X,Y\subseteq V(H)$ satisfy $|X|n_1=|Y|n_2$ and either $X\sqcup Y=V(H)$ or $X=Y=V(H)$. Assume every $\sim_H^Y$-class in $X$ has size $1$ and let $\e>0$ satisfy
    $$
    \e<\left(\frac{1}{2(k^2+1)|Y|}\right)^{2^{k}}.
    $$
  Then any $\e$-regular partition $\calP$ of any $(X,Y;n_1,n_2)$-blowup $\mathbf{H}$ of $H$ satisfies
   $$
   |\calP|\geq  (1-2(k^2+1)\e^{1/2^k})(1- \e^{1/4})|X|.
   $$
  
\end{corollary}
\begin{proof}
Say ${\bf H}$ has vertex set $\bigcup_{u\in V(H)}Z_u$, and set ${\bf X}=\bigcup_{u\in X}Z_u$.  Let $\calP'$ be as in the conclusion of Lemma \ref{lem:blowup}. Let $X_1,\ldots, X_t$ enumerate the $\sim_H^Y$ classes in $X$, and for each $1\leq i\leq t$, set ${\bf X}_i=\bigcup_{u\in X_i}Z_u$. Conclusion (ii) of Lemma \ref{lem:blowup} implies that for all $P\in \calP'$, there is $1\leq f(P)\leq t$ such that $|P\cap {\bf X}_{f(P)}|\geq (1- \e^{1/4})|P\cap {\bf X}|$.  From this we can deduce that for all $P\in \calP'$,
$$
(1- \e^{1/4})|P\cap {\bf X}|\leq  |P\cap {\bf X}_{f(P)}|\leq |{\bf X}_{f(P)}|=n_1,
$$
where the final equality is by definition of an $(X,Y; n_1,n_2)$-blowup, and our assumption that all $\sim_H^Y$-classes have size $1$.  On the other hand, conclusion (i) of Lemma \ref{lem:blowup} implies $|\bigcup_{P\in \calP'}P\cap {\bf X}|\geq (1-2(k^2+1)\e^{1/2^k})|{\bf X}|$.   Combining these observations, we deduce the following.
$$
|\calP|\geq \frac{|\bigcup_{P\in \calP'}P\cap {\bf X}|}{\max_{P\in \calP'}|P\cap {\bf X}|}\geq \frac{(1-2(k^2+1)\e^{1/2^k})|{\bf X}|}{(1- \e^{1/4})^{-1}n_1}=  (1-2(k^2+1)\e^{1/2^k})(1- \e^{1/4})|X|,
$$
where the last equality is by definition of ${\bf X}$ and the definition of an $(X,Y; n_1,n_2)$-blowup.
\end{proof}

\section{Jump between polynomial and exponential}

In this section we characterize the jump between polynomial and exponential growth.  

\begin{proposition}\label{prop:polyexp}
Suppose $\calH$ is a hereditary $k$-graph property. Then one of the following holds.
\begin{enumerate}
\item $\calH$ is close to finite VC-dimension and $M_{\calH}(\e)\leq \e^{-C}$ for some $C>0$, or
\item $\calH$ is far from finite VC-dimension and $M_{\calH}(\e)\geq 2^{\e^{-C}}$ for some $C>0$.
\end{enumerate}
\end{proposition}
\begin{proof}
First suppose $\cH$ is close to finite VC-dimension. Then there is $\calH'$ with finite VC-dimension such that $\cH$ is close to $\cH'$. By Theorem \ref{thm:FPS} and Proposition \ref{prop:homogimpliesreg}, there are constants $c, c'$ such that $M_{\cH'}(\e)\leq c\e^{-c'}$. By Fact \ref{fact:close}, this implies there is a constant $C$ such that $M_{\cH}(\e)\leq \e^{-C}$. 

Suppose now $\cH$ is far from finite VC-dimension. Let $\e>0$ be sufficiently small compared to $k$, and set $K=\lceil (2(k^2+1) (k-1)\e^{1/2^{k}})^{-1}\rceil-1$.  By Proposition \ref{prop:vcclose}, there exists some $H\in {\bf B}_{\calH}\cap \calU_k(K)$. By definition of $\calU_k(K)$, $V(H)$ has the form $X\sqcup Y$, where 
$$
X=\{c_S:S\subseteq [K]\}\text{ and }Y=\{a_i^j: 1\leq i\leq k-1, 1\leq j\leq K\},
$$
and edge set $E(H)$ satisfying that for all $j\in [K]$ and $S\subseteq [K]$,
\begin{align*}
\text{ $j\in S\Rightarrow a_1^{j}\ldots a_{k-1}^jc_S\in E(H)$ and $j\notin S\Rightarrow a_1^{j}\ldots a_{k-1}^jc_S\notin E(H)$.}
\end{align*}
Clearly $|Y|\leq (k-1)K$, while the above condition on $E(H)$ implies $|X|=2^K$.  Now let $n_1$ be sufficiently large and divisible by $|Y|$ and define $n_2=n_1|X||Y|^{-1}$.  We note that our assumptions on the edge set of $H$ imply that every $\sim_H^Y$-class in $X$ has size equal to $1$.  Since $H\in {\bf B}_{\calH}$ and $\calH$ is hereditary, it is not hard to see there must exist some $(X,Y; n_1,n_2)$-blowup ${\bf H}$ of $H$ in $\calH$.  

We now check that the hypotheses of Corollary \ref{cor:blowup} hold.  By construction, $X\sqcup Y=V(H)$, $n_1|X|=n_2|Y|$, and every $\sim_H^Y$-class in $X$ has size $1$ by above. By the definition of $K$ and since $|Y| \leq (k-1)K$, we have
$$
\e<\left(\frac{1}{2(k^2+1)|Y|}\right)^{2^{k}}.
$$
Consequently, Corollary \ref{cor:blowup} implies any $\e$-regular partition $\calP$ of ${\bf H}$ satisfies 
$$
|\calP|\geq  (1-2 (k^2+1)\e^{1/2^k}) (1- \e^{1/4})|X|=(1-2 (k^2+1)\e^{1/2^k})(1- \e^{1/4})2^K.
$$
Recalling the definition of $K$, we see this implies that for some constant $C$ depending only on $k$, $|\calP|\geq 2^{\e^{-C}}$.  This completes the proof.
\end{proof}

\section{Constant to Polynomial Jump}

In this section we characterize the jump between constant and polynomial.  We will be using throughout a tool introduced in Section \ref{sec:LB}, namely the equivalence relation on the vertices of a $k$-graph obtained by grouping together twins (see Definition \ref{def:twins}).   In this section, $k$-graphs where these equivalence classes are especially small will play a crucial role.   For our purposes, the central notion is where we allow the $\sim_H$-classes to have size at most $k$. We will call this notion \emph{almost primeness} to distinguish it from pre-existing notions of prime hypergraphs. 

\begin{definition}
    Suppose $H = (V, E)$ is a $k$-graph. Then $H$ is \emph{almost prime} if every $\sim_H$-class has size at most $k$.
\end{definition}

 The main result of this section is the following.
 
 \begin{proposition}\label{prop:polyconstant}
 Suppose $\calH$ is a hereditary $k$-graph property. Then one of the following holds.
 \begin{enumerate}
 \item ${\bf B}_{\calH}$ contains infinitely many non-isomorphic almost prime $k$-graphs. In this case, $M_{\calH}(\e)\geq \e^{-C}$ for some constant $C>0$.
 \item ${\bf B}_{\calH}$ contains finitely many non-isomorphic almost prime $k$-graphs. In this case,  $M_{\calH}(\e)=C$, where $C$ is the maximum number of $\sim_H$-classes appearing in an almost prime element $H$ of ${\bf B}_{\calH}$.
 \end{enumerate}
 \end{proposition}

 We will prove part (1) in Subsection \ref{ss:polyLB}. In Subsection \ref{ss:twins}, we will prove several lemmas about almost prime $k$-graphs and twin classes, with the aim of showing one can characterize any fixed upper bound on the number of twin classes by finitely many forbidden $k$-graphs. This is required to obtain the finite determination of growth classes described after the statement of Theorem \ref{thm:main}.  In  Subsection \ref{ss:constantUB} we prove part (2) of Proposition \ref{prop:polyconstant}, after which we will prove our main theorem, Theorem \ref{thm:main}.
 
\subsection{Polynomial lower bound}\label{ss:polyLB}
Our goal in this subsection is to prove part (1) of Proposition \ref{prop:polyconstant}.  We begin by setting notation for a special collection of $k$-graphs.   

\begin{definition}
    Let $\ell\geq 1$. Define $\calA(\ell)$ to be the class of $k$-graphs $H$ with vertex set of the form 
    $$
    V (H)=\{a_i^j: 1\leq j\leq k, 1\leq i\leq \ell\},
    $$
   satisfying
   \begin{enumerate}[(i)]
   \item  $|\{a_i^1,\ldots, a_i^k\}|=k$ for each $1\leq i\leq \ell$, and 
   \item $\{a_i^1: 1\leq i\leq \ell\}\cap \{a_i^j: 1\leq i\leq\ell, 2\leq j\leq k\}=\emptyset$,
\end{enumerate}
    and with edge set $E$ satisfying one of the following. 
    \begin{enumerate}[(a)]
        \item  $a_i^1a_j^2\ldots a_j^k \in E$ if and only if $i=j$
        \item  $a_i^1a_j^2\ldots a_j^k \in E$ if and only if $i\neq j$
        \item  $a_i^1a_j^2\ldots a_j^k \in E$ if and only if $i\leq j$.
    \end{enumerate}
\end{definition}

The $k$-graphs in $\calA(\ell)$ will play an important role in this section for two reasons. First,  their blowups will suffice for the desired lower bound. Secondly, one can always find elements of $\calA(\ell)$ in sufficiently large almost prime $k$-graphs, as the following lemma states.

\begin{lemma}\label{lem:A(l)}
For all $\ell\geq 1$ there is $N\geq 1$ such that any almost prime $k$-graph with at least $N$ vertices contains an induced sub-$k$-graph which is an element of $\calA(\ell)$.
\end{lemma}

The proof of Lemma \ref{lem:A(l)} is an adaptation of the 3-uniform case from \cite{part2} and appears in the appendix.  Despite the fact that Lemma \ref{lem:A(l)} connects $\calA(\ell)$ to almost prime $k$-graphs in an important way, it is not true that all elements of $\calA(\ell)$ are themselves almost prime. Indeed, assume $\ell>k\geq 3$, and consider the $k$-graph $H$ with vertex set $V(H)=\{a_i^j: 1\leq j\leq k, 1\leq i\leq \ell\}$ and edge set $E(H)= \{a_i^1\ldots a_i^{k-1}a_j^k: 1\leq i,j\leq \ell\} $. By definition $H\in \calA(\ell)$, however $H$ is not almost prime, as the elements in $a^k_1,\ldots, a^k_{\ell}$ are in a single $\sim_H$-class.

We are now ready to prove what we need for part (1) of Proposition \ref{prop:polyconstant}.

 \begin{theorem}
    Let $\cH$ be a hereditary $k$-graph property such that $\bf B_\calH$ contains infinitely many non-isomorphic almost prime $k$-graphs. Then $M_{\calH}(\e)\geq \Omega(\e^{-1/2^k})$.
\end{theorem}

\begin{proof}
    Let $\e>0$ be sufficiently small and set
    $$
    \ell = \left\lceil\left(2(k^2+1)(k-1)\e^{1/2^k}\right)^{-1}\right\rceil-1.
    $$ 
    Since ${\bf B}_{\calH}$ contains infinitely many non-isomorphic almost prime $k$-graphs, Lemma \ref{lem:A(l)} implies there is $H\in \bf B_\calH$ such that $H$ contains an element of $\calA(\ell)$ as an induced sub-$k$-graph.  Since ${\bf B}_{\calH}$ is hereditary (see Fact \ref{fact:bhhgraphs}), this implies there is some $A\in \calA(\ell)\cap \bf B_{\calH}$. We wish to apply Corollary \ref{cor:blowup} to the $k$-graph $A$ and appropriately chosen sets $X$ and $Y$. By definition of $\calA(\ell)$, $A$ has vertex set of the form 
    $$
    \{a_i^j: 1\leq j\leq k, 1\leq i\leq \ell\},
    $$
    and edge set $E$ satisfying one of the following.
    \begin{enumerate}
        \item  $a_i^1a_j^2\ldots a_j^k \in E$ if and only if $i=j$
        \item  $a_i^1a_j^2\ldots a_j^k \in E$ if and only if $i\neq j$
        \item  $a_i^1a_j^2\ldots a_j^k \in E$ if and only if $i\leq j$
    \end{enumerate}
    Let $X = \{a_1^1, \ldots a_\ell^1\}$ and let $Y = \{a_i^j:i\in[\ell], j\in\{2, \ldots, k\}\}$. Note $|Y|\leq (k-1)\ell$ by definition, and the conditions on $E$ imply $|X|=\ell$. By our choice of $\ell$, we have
    $$
    \e<\left(\frac{1}{2(k^2+1)(k-1)\ell}\right)^{2^k}\leq \left(\frac{1}{2(k^2+1)|Y|}\right)^{2^k}.
    $$
    Let $n_1, n_2$ be sufficiently large integers satisfying $n_1|X| = n_2|Y|$. Since $A\in \bf B_{\calH}$, there is an  $(X, Y; n_1, n_2)$-blowup $\bf A$ of $A$ in $\calH$. Say $\bf A$ has vertex set $V(\mathbf{A}) = \cup_{u\in V(A)}Z_u$ where $|Z_u| = n_1$ for each $u\in X$ and $|Z_u| = n_2$ for each $u\in Y$. Since $A\in \calA(\ell)$, all of the $\sim_A^Y$ classes in $X$ have size $1$. Thus we have satisfied the hypotheses of Corollary \ref{cor:blowup}. Consequently,  any $\e$-regular partition $\calP$ of $\bf A$ has size at least 
    \begin{align*}
        |\calP|\geq  (1-2(k^2+1)\e^{1/2^k}) (1- \e^{1/4})|X|&=(1-2(k^2+1)\e^{1/2^k}) (1- \e^{1/4})\ell.   
    \end{align*}
   Recalling the definition of $\ell$, we see this implies $M_\calH(\e)\geq\Omega(\e^{-1/2^k})$ as desired.
\end{proof}

\subsection{Bounds on twin classes}\label{ss:twins} The goal of this subsection is to prove that we can characterize the class of $k$-graphs $H$ admitting a fixed bound $C$ on the number of $\sim_H$-classes via finitely many forbidden $k$-graphs.  As this will require a bit of work, we first explain the purpose of such a result.   Recalling that $\calG^{(k)}$ denotes the class of all finite $k$-graphs, we define the following class for each integer $t>0$:  
$$
\calH_t^{(k)}=\{H\in \calG^{(k)}:\text{ $H$ has at most $t$-many $\sim_H$ classes}\}. 
$$

It is not hard to see $\calH_t^{(k)}$ is a hereditary $k$-graph property. Indeed, deleting vertices can never increase the number of twin classes (see Lemma \ref{lem:coarsening} below).  Explicitly characterizing $\calH_t^{(k)}$ by an infinite family of forbidden $k$-graphs is also very easy.  Specifically, if we set 
$$
\calS_t^{(k)}=\{H\in \calG^{(k)}: \text{ $H$ has at least $t$-many $\sim_H$-classes}\},
$$
then the following holds by definition.
\begin{fact}\label{fact:st}
For all integers $t>0$, $\calH_t^{(k)}=\Forb(\calS_{t+1}^{(k)})$.
\end{fact}
In Subsection \ref{ss:constantUB}, we will see that the class $\calH_t^{(k)}$ plays a crucial role in the constant range of our main theorem, Theorem \ref{thm:main}.  Inspecting that proof, one can see the same statement appearing in Theorem \ref{thm:main} can be proved using Fact \ref{fact:st}.  This boils down to the fact that Fact \ref{fact:blowup} applies to both infinite and finite sets of forbidden $k$-graphs.  However, we are interested in proving a stronger version of Theorem \ref{thm:main} in which the upper bounds arise from explicitly defined  \emph{finite} families of forbidden $k$-graphs.  For this, we require a characterization of $\calH_t^{(k)}$ by \emph{finitely many} forbidden $k$-graphs, and this is the goal of the current subsection.  Specifically, we prove Lemma \ref{lem:forb} below, where for an integer $t>0$, 
\begin{align*}
\calF_t^{(k)}=\{H\in \calG^{(k)}:\text{$H$ is almost prime, $H$ has at least $t$-many $\sim_H$-classes, }&\\
 \text{ and $H$ has at most $k^2t$ vertices}\}&.
\end{align*}

\begin{lemma}\label{lem:forb}
For all integers $t>0$, $\calH_t^{(k)} = \Forb(\calF_{t+1}^{(k)})$. 
\end{lemma}

Our proof of Lemma \ref{lem:forb} will follow closely the strategy used in the $3$-uniform setting in \cite{part2}.   As such we will require several preliminaries about twin classes.  First, we observe that by  definition, for any almost prime $k$-graph $H$  with exactly $\ell$ many $\sim_H$-classes, one has $|V(H)|\leq k\ell$.  Our next lemma follows immediately from the definition of $\sim_H$ (see Definition \ref{def:twins}).

\begin{lemma}\label{lem:coarsening}
    Let $H = (V, E)$ be a $k$-graph, and let $U_1, \ldots, U_\ell$ be the $\sim_H$-classes of $H$. Let $H' = (V', E')$ be an induced sub-$k$-graph of $H$. Then $\sim_{H'}$-classes of $H'$ are unions of sets from $\{U_i\cap V':i\in[\ell]\}$.
\end{lemma}

The following lemma tells us we can witness the number of $\sim_H$-classes appearing in a $k$-graph $H$ using an almost prime induced sub-$k$-graph.  This is the $k$-ary analogue of Lemma B.8  in \cite{part2}.

\begin{lemma}\label{lem:representative}
    Let $H = (V, E)$ be a $k$-graph and suppose $H$ has exactly $\ell$ many $\sim_H$-classes.  Then there is an induced sub-$k$-graph $H' = (V', E')$ of $H$ which is almost prime and also has exactly $\ell$ many $\sim_{H'}$-classes.
\end{lemma}
\begin{proof}
    Let $U_1, \ldots, U_\ell$ be the $\sim_H$-classes of $H$. For each $i\in[\ell]$ choose any $U_i'\subseteq U_i$ such that $|U_i'| = \min\{|U_i|, k\}$ and let 
    $$
    H' = H\left[\bigcup_{i=1}^\ell U_i'\right].
    $$
    It suffices to show the $\sim_{H'}$-classes are exactly $U_1', \ldots, U_\ell'$. By Lemma \ref{lem:coarsening}, each $U_i'$ is contained in some $\sim_{H'}$-class. Thus it only remains to show that for any $i\neq j$, $U_i'$ and $U_j'$ are not contained in the same $\sim_{H'}$-class. 

    Fix $1\leq i\neq j\leq \ell$. Let $u\in U_i'$ and $v\in U_j'$.  Our goal is to show $u\nsim_{H'}v$. Since $u\not\sim_H v$, there are $w_1, \ldots, w_{k-1}\in V\setminus \{u,v\}$ such that (after possibly relabeling), $w_1\ldots w_{k-1}u\in E$ and $w_1 \ldots w_{k-1}v\notin E$.  
    
    For each $\alpha\in [\ell]$, let $S_\alpha = \{j\in[k-1]:w_j\in U_\alpha\}$.  We claim $|U'_{\alpha}\setminus \{u,v\}|\geq |S_{\alpha}|$.  Note that $\{w_1,\ldots,w_{k-1}\}\cap \{u,v\}=\emptyset$ implies 
    \begin{align}\label{al:ui1}
    |U_{\alpha}\setminus \{u,v\}|\geq |U_{\alpha}\cap \{w_1,\ldots, w_{k-1}\}|=|S_{\alpha}|.
     \end{align}
    If $|U_{\alpha}'|=|U_{\alpha}|$, then $U_{\alpha}'=U_{\alpha}$, and (\ref{al:ui1}) immediately implies $|U_{\alpha}'\setminus \{u,v\}|\geq |S_{\alpha}|$.  If $|U_{\alpha}'|\neq |U_{\alpha}|$, then $|U_{\alpha}'|$ must be equal to $k$. In this case, we observe that at most one of $u$ or $v$ can be contained in $U_\alpha'$, and consequently,
   \begin{align}\label{al:ui}
   |U_\alpha'\setminus \{u,v\}|\geq |U_\alpha'|-1=k-1\geq |S_{\alpha}|,
   \end{align}
  where the last inequality is because $S_{\alpha}\subseteq [k-1]$ by definition.  This concludes our verification that for all $\alpha\in [\ell]$, $|U'_{\alpha}\setminus \{u,v\}|\geq |S_{\alpha}|$.   Thus we may choose distinct vertices $\{x_j:j\in S_\alpha\}$ from $U_\alpha'\setminus\{u,v\}$, letting $x_j=w_j$ in cases where $w_j\in U_\alpha'$.  Since $x_1,\ldots, x_{k-1}\in V(H')$, we will be done after we prove the following claim (since it will imply $u\nsim_{H'}v$). 

    \begin{claim}
        We have $x_1\ldots x_{k-1}u\in E$ and $x_{1}\ldots x_{k-1}v\notin E$. 
    \end{claim}
    \begin{proof}
Set $x(1)=\emptyset$, and for each $1<j\leq k$, set $x(j)=\{x_1,\ldots, x_{j-1}\}$.  Then define $w(k)=\emptyset$, and for each $1\leq j<k$, set $w(j)=\{w_j,\ldots, w_{k-1}\}$.  Note that in this notation, $x(1)\cup w(1)=\{w_1,\ldots, w_{k-1}\}$ and $x(k)\cup w(k)=\{x_1,\ldots, x_{k-1}\}$.  Our assumptions tell us 
    \begin{align}\label{al:xwu}
 x(1)\cup w(1)\cup \{u\}\in E\text{ and }x(1)\cup w(1)\cup \{v\}\notin E.
    \end{align}
For all $1\leq j\leq k-1$, we have $w_j\sim_H x_j$, and consequently, 
   $$
   \{x_1,\ldots, x_{j-1},x_j,w_{j+1}, \ldots, w_{k-1},v\}\in E \Leftrightarrow \{x_1,\ldots, x_{j-1},w_j,w_{j+1},\ldots, w_{k-1},v\}\in E.
   $$
   In other words, $x(j+1)\cup w(j+1)\cup \{v\}\in E\Leftrightarrow x(j)\cup w(j)\cup \{v\}\in E$.  A similar argument shows that for all $1\leq j\leq k-1$, $x(j+1)\cup w(j+1)\cup \{u\}\in E\Leftrightarrow x(j)\cup w(j)\cup \{u\}\in E$.  Combining with the induction hypothesis, we have $x(j+1)\cup w(j+1)\cup \{u\}\in E$ and $x(j+1)\cup w(j+1)\cup \{v\}\notin E$.  Combining these observations with (\ref{al:xwu}), we have that   $x(k)\cup w(k)\cup \{u\}\in E$ and $x(k)\cup w(k)\cup \{v\}\notin E$.  Translating, this means  $x_1\ldots x_{k-1}u\in E$ and $x_{1}\ldots x_{k-1}v\notin E$, as desired. 
    \end{proof}
\end{proof}

Next we have a deletion lemma, which says that, given a $k$-graph $H$, we can delete a $\sim_H$-class without collapsing too many other $\sim_H$ classes. A weaker version of this was proved in the $3$-uniform setting in \cite{part2} (see Lemma B.10 there).  While that proof can likely be adapted to the $k$-uniform setting, we have opted to present a shorter proof of a stronger statement generated by ChatGPT (rewritten by us).

\begin{lemma}\label{lem:deletion}
    Let $H = (V, E)$ be a $k$-graph, let $\ell\geq 2$ be an integer, let $U_1, \ldots, U_\ell$ enumerate the $\sim_H$-classes of $H$, and for each $i\in[\ell]$, let $H_i := H[V\setminus U_i]$. Then there is $i\in [\ell]$ such that $H_i$ has at least $\frac{\ell-1}{k}$-many $\sim_{H_i}$-classes.
\end{lemma}

\begin{proof}
    For each $i\in[\ell]$, fix $u_i\in U_i$. For each $1\leq s<t\leq \ell$, since $u_s\nsim_H u_t$, we can fix a set $W_{st} = \{w^{st}_1, \ldots, w^{st}_{k-1}\}\subseteq V\setminus\{u_s,u_t\}$ such that $w^{st}_1 \ldots w^{st}_{k-1}u_s\in E$ and $w^{st}_1 \ldots w^{st}_{k-1}u_t\notin E$, or vice versa.

For each $i\in [\ell]$, define $\calP_i$ to be the partition of $[\ell]\setminus\{i\}$ such that $s, t$ are in the same part of $\calP_i$ if and only if $U_s, U_t$ lie in the same $\sim_{H_i}$-class. Setting  $\ell_i = |\calP_i|$, we see that $\ell_i$ is the number of $\sim_{H_i}$-classes in $H_i$. Further define 
    \begin{align}\label{al:midef}
    M_i = \sum_{P\in \calP_i}\binom{|P|}{2}.
    \end{align}
    Then $M_i$ is the total number of pairs of $\sim_H$-classes which merged in $\sim_{H_i}$.  Next, for each $1\leq s<t\leq \ell$,  define 
    $$
    I_{st}:=\{j\in [\ell]\setminus\{s, t\}: \text{ $s, t$ lie in the same part of $\calP_j$}\}.
    $$
     We claim if $i\in I_{st}$, then $W_{st}\cap U_i\neq \emptyset$. Indeed, suppose $W_{st}\cap U_i = \emptyset$. Then $w_1^{st}, \ldots, w_{k-1}^{st}\in V(H_i)$, and therefore $W_{st}$ separates $U_s$ and $U_t$ in $H_i$. Thus $i\notin I_{st}$. Consequently, since each $W_{st}$ contains only $k-1$ elements, it intersects at most $k-1$ of the sets $U_i$. Thus, for each $s<t\in [\ell]$, $|I_{st}|\leq k-1$. Note next that $M_i = |\{(s, t):1\leq s<t\leq \ell, i\in I_{st}\}|$. Therefore, we have
    $$
    \sum_{i\in[\ell]}M_i = \sum_{1\leq s<t\leq \ell}|I_{st}|\leq \binom{\ell}{2}(k-1).
    $$
    By the pigeonhole principle, there is some $i\in[\ell]$ such that 
    \begin{align}\label{al:mi}
    M_i\leq \frac{1}{\ell}\binom{\ell}{2}(k-1) = \frac{(\ell-1)(k-1)}{2}.
    \end{align}
    Enumerate the parts of the partition $\calP_i$ by $P_1, \ldots, P_{\ell_i}$. Then, by definition of $\calP_i$,
    \begin{align}\label{al:pi}
    \sum_{j\in[\ell_i]}|P_j| = \ell-1.
    \end{align}
    We also have $|P_j|^2 = 2\binom{|P_j|}{2}+|P_j|$ and therefore, combining  (\ref{al:mi}), (\ref{al:pi}), and the definition of $M_i$ (see (\ref{al:midef})), we have 
    $$
    \sum_{j\in[\ell_i]}|P_j|^2= 2M_i + \ell-1\leq (k-1)(\ell-1)+\ell-1=k(\ell-1).
    $$
 Consequently,
    $$
    \ell_i\sum_{j\in[\ell_i]}|P_j|^2\leq \ell_ik(\ell-1).
    $$
    On the other hand, by Cauchy-Schwarz,
    $$
    \ell_i\sum_{j=1}^{\ell_i}|P_j|^2\geq (\ell-1)^2.
    $$
    Putting these two inequalities together, we have $\ell_ik(\ell-1)\geq (\ell-1)^2$, so $\ell_i\geq \frac{\ell-1}{k}$.    This says exactly that the $k$-graph $H_i$ has at least $\frac{\ell-1}{k}$-many $\sim_{H_i}$ classes.
\end{proof}

Our next lemma says that if $H$ is a $k$-graph with at least $Ck$-many $\sim_H$ classes, then we can always find an induced almost prime sub-$k$-graph $H'$ with at least $C$-many $\sim_{H'}$-classes and a bounded number of vertices.  This is a straightforward adaptation of the $3$-uniform analogue (see Lemma B.12 in \cite{part2}). 

\begin{lemma}\label{lem:bounded}
    Let $C>0$ be an integer. Let $H$ be a $k$-graph and suppose $H$ has at least $Ck$-many $\sim_H$ classes. Then there is an induced sub-$k$-graph $H'$ of $H$ such that $H'$ is almost prime, has at least $C$-many $\sim_{H'}$ classes, and has at most $k^2C$-many vertices.
\end{lemma}

\begin{proof}
    We will repeatedly apply Lemma \ref{lem:deletion}.  Let $\ell$ be the number of $\sim_H$-classes in $H$.  We construct an integer $t\geq 0$ and a sequence $H_0,\ldots, H_t$ with the goal of producing $H_t$ having at least $C$, and at most $kC$-many $\sim_{H_t}$-classes.   
    
    Step $0$: Let $H_0 = H$. By assumption, $H_0$ has at least $Ck$-many $\sim_{H_0}$-classes. If $H_0$ also has at most $Ck$-many $\sim_{H_0}$-classes, set $t=0$ and end the construction. Otherwise, $H_0$ has at least $(Ck+1)$-many $\sim_{H_0}$-classes.  Go to the next step. 
    
  Step $i+1$:  Assume $i\geq 0$, and suppose by induction we have defined a $k$-graph $H_i$ with at least $(Ck+1)$-many, and at most $(\ell-i)$-many $\sim_{H_i}$ classes. Apply Lemma \ref{lem:deletion} to $H_i$ to obtain an induced sub-$k$-graph $H_{i+1}$ with at least $C$-many, and at most $(\ell-i-1)$-many $\sim_{H_{i+1}}$-classes. If $H_{i+1}$ has at most $Ck$-many $\sim_{H_{i+1}}$-classes, set $t=i+1$ and end the construction. Otherwise, $H_{i+1}$ has at least $(Ck+1)$-many $\sim_{H_{i+1}}$-classes, and we go to the next step.  
    
Since the number of $\sim$-classes decreases at each step in this construction, it will halt at some finite stage $t$, yielding an induced sub-$k$-graph $H_t$ of $H$ with at least $C$-many, and at most $Ck$-many $\sim_{H_t}$-classes. Apply Lemma \ref{lem:representative} to the $k$-graph $H_t$ to obtain an almost prime induced sub-$k$-graph $H_t'$ of $H_t$ with at least $C$-many, and at most $Ck$-many $\sim_{H_t'}$-classes.  Since $H_t'$ is almost prime and has at most $kC$-many $\sim_{H_t'}$ classes, it has at most $k^2C$-many vertices.
\end{proof}

 We are now ready to prove Lemma \ref{lem:forb}.
 
 \vspace{2mm}
 
 \noindent{\bf Proof of Lemma \ref{lem:forb}.} First note that if $F\in \calF_{t+1}^{(k)}$, then $F$ has at least $(t+1)$-many $\sim_F$ classes,  so $F\notin \calH_t^{(k)}$. This shows $\calH_t^{(k)}\subseteq \Forb(\calF_{t+1}^{(k)})$.  
    
 For the reverse containment, suppose $H\notin \calH_t^{(k)}$.  Our goal is to show $H\notin \Forb(\calF_{t+1}^{(k)})$, i.e. that $H$ contains an induced sub-$k$-graph in $\calF_{t+1}^{(k)}$.  Let $\ell$ be the number of $\sim_H$-classes in $H$. Since $H\notin \calH_t^{(k)}$, $\ell\geq t+1$.  By Lemma \ref{lem:representative}, there is an induced almost prime sub-$k$-graph $H'$ of $H$ with $\ell$-many $\sim_{H'}$-classes. Then $|V(H')|\leq k\ell$. If it is also the case that $\ell\leq k(t+1)$, then $H'$ has at least $(t+1)$-many $\sim_{H'}$-classes  and  at most $k^2(t+1)$-vertices. Thus in this case, we have $H'\in \calF_{t+1}^{(k)}$ by definition, so we are done.
 
Assume now that $\ell>k(t+1)$. Then Lemma \ref{lem:bounded} implies there is an induced sub-$k$-graph $H''$ of $H'$ such that $H''$ has at least $(t+1)$-many $\sim_{H''}$-classes and at most $k^2(t+1)$ vertices. In this case we have $H''\in \calF_{t+1}^{(k)}$, so we are done. 
\qed

\subsection{Constant bound}\label{ss:constantUB} In this subsection we prove what we need for part (2) of Proposition \ref{prop:polyconstant}, after which we will prove the main theorem, Theorem \ref{thm:main}. 

As far as Proposition \ref{prop:polyconstant} is concerned, our goal is to show that when $\bf B_\calH$ contains only finitely many non-isomorphic almost prime $k$-graphs, then there is a constant $C$ such that $M_{\cH}(\e)=C$.   We begin by proving that if $\bf B_\calH$ contains only finitely many non-isomorphic almost prime $k$-graphs, then $M_\calH(\e)\leq C$, where $C$ is the maximum number of twin classes appearing in any almost prime element of ${\bf B}_{\calH}$.  

\begin{theorem}\label{thm:constUB}
    Let $\calH$ be a hereditary $k$-graph property. Suppose that $\bf B_{\calH}$ contains finitely many non-isomorphic almost prime $k$-graphs. Let $C\geq 1$ be the maximum value such that $\bf B_\calH$ contains an almost prime $k$-graph with $C$-many $\sim$-classes. Then for all sufficiently small $\e>0$, $M_\calH(\e) \leq C$.
\end{theorem}
\begin{proof}
Fix $\e>0$ sufficiently small compared to $C$ and $k$.   By Lemma \ref{lem:forb}, $\calH_C^{(k)} = \Forb(\calF_{C+1}^{(k)})$, and by assumption, ${\bf B}_{\calH}\cap \calF_{C+1}^{(k)}=\emptyset$. By Fact \ref{fact:blowup}, this implies $\calH$ is close to $\calH_C^{(k)} $. Therefore, it suffices to show $M_{\calH_C^{(k)}}(\e)\leq C$ by Fact \ref{fact:close}.

    Let $H$ be a sufficiently large element of $\calH_C^{(k)}$.  We show $H$ has an $\e$-regular partition with at most $C$ parts. By definition, $H$ has at most $C$-many $\sim_H$-classes. Let $\calU$ be the set of $\sim_H$-classes of $H$, and note $|\calU|\leq C$. We show $\calU$ is $\e$-regular with respect to $H$. Let 
    $$
    \Sigma_{\textrm{small}}=\{(U_1,\ldots, U_k)\in \calU^k: \text{ for some }i\in [k], |U_i|\leq \e^2|V(H)|\}.
    $$
   Then 
    $$
 \left|\bigcup_{(U_1,\ldots, U_k)\in \Sigma_{\textrm{small}}}U_1\times \cdots \times U_k\right|\leq \e^2kC|V(H)|^k\leq \e |V(H)|^k,
    $$
where the last inequality is because $\e$ is sufficiently small compared to $k$ and $C$.   It now suffices to show every $(U_1,\ldots, U_k)\in \calU^k\setminus  \Sigma_{\textrm{small}}$ is $\e$-regular with respect to $H$.  Fix $(U_{1}, \ldots, U_{k})\in \calU^k\setminus  \Sigma_{\textrm{small}}$. Since each of the sets $U_1,\ldots, U_k$ are $\sim_H$-classes, it is an exercise to prove that one of the following holds: 
 \begin{enumerate}
 \item $ |\overline{E(H)}\cap (U_{1}\times \ldots \times U_{k})|=0$ or 
 \item $|\overline{E(H)}\cap (U_{1}\times \ldots \times U_{k})|= |\{(u_1,\ldots, u_k)\in U_1\times \cdots \times U_k: |\{u_1,\ldots, u_k\}|=k\}|.$
  \end{enumerate}
Clearly in case (1), $(U_1,\ldots, U_k)$ is $\e$-regular with respect to $H$. So let us assume we are in case (2).  In this case, an elementary computation shows 
$$
d_H(U_1,\ldots, U_k)\geq ( 1-k2^k\max\{|U_i|^{-1}: i\in [k]\})\geq 1-k 2^k \e^{-2}|V(H)|^{-1},
$$
where the last inequality uses that $(U_1,\ldots, U_k)\notin \Sigma_{\textrm{small}}$.  Since we assumed $H$ is sufficiently large, this implies $d_H(U_1,\ldots, U_k)\geq 1-\e^{2^{k+1}}$, and consequently, $(U_1,\ldots, U_k)$ is $\e$-regular with respect to $H$ by Proposition \ref{prop:homogimpliesreg}.  We have now shown every sufficiently large element of $\calH_C^{(k)}$ has an $\e$-regular partition with at most $C$ parts, and thus $M_{\calH_C^{(k)}}(\e)\leq C$, as desired.
\end{proof}

Next, we show that when  $\bf B_{\calH}$ contains finitely many non-isomorphic almost prime $k$-graphs, $M_{\calH}$ is equal to the constant appearing in Theorem \ref{thm:constUB}. 

\begin{theorem}\label{thm:constmain}
    Let $\calH$ be a hereditary $k$-graph property. Suppose that $\bf B_{\calH}$ contains finitely many non-isomorphic almost prime $k$-graphs. Let $C\geq 1$ be the maximum value such that $\bf B_\calH$ contains an almost prime $k$-graph with $C$-many $\sim$-classes. Then for all sufficiently small $\e>0$, $M_\calH(\e) =C$.
\end{theorem}
\begin{proof}
Fix $\e>0$ sufficiently small compared to $C$ and $k$.  That $M_{\calH}(\e)\leq C$ follows immediately from Theorem \ref{thm:constUB}.

 Turning to the lower bound, let $H\in \bf B_{\calH}$ be an almost prime $k$-graph with $C$-many $\sim_H$-classes.  Enumerate the $\sim_H$ classes by $U_1, \ldots, U_C$. Since $H$ is almost prime, $|V(H)|\leq kC$. Let $n$ be any integer which is sufficiently large compared to $\e^{-1}$. Since $H\in \bf B_{\calH}$, there is an $n$-blowup $\bf H$ of $H$ with ${\bf H}\in \calH$. Say $\bf H$ has vertex set 
    $V(\mathbf{H}) = \bigcup_{u\in V(H)}Z_u$
    where $|Z_u| = n$ for all $u\in V(H)$. 
    
    We wish to apply Lemma \ref{lem:blowup} to $\bf H$. Set $X = Y = V(H)$ and $n_1 = n_2 = n$. Clearly $n_1|X| = n_2|Y|$, and  $Cn\leq |V(\mathbf{H})|\leq kCn$. Since $\e$ is sufficiently small, we may assume 
    $$
    \e< \left(\frac{1}{2(k^2+1)kC}\right)^{2^k}\leq \left(\frac{1}{2(k^2+1)|V(H)|}\right)^{2^k}.
    $$
    Thus the hypotheses of Lemma \ref{lem:blowup} are satisfied. Let $\calP$ be an $\e$-regular partition of $\bf H$. Then by Lemma \ref{lem:blowup}, there is $\calP'\subseteq \calP$ such that the following hold, where for each $i\in [C]$ let $\mathbf{U}_i = \bigcup_{u\in U_i}Z_u$. 
    \begin{enumerate}[(i)]
        \item $\left|\bigcup_{P\in \calP'} P\right|\geq (1-2(k^2+1)\e^{1/2^k})|V(\mathbf{H})|$
    and
    \item For all $P\in \calP'$, there is $1\leq \sigma(P)\leq C$ such that $|P\cap \mathbf{U}_{\sigma(P)}|\geq (1-\e^{1/4})|P|$. 
    \end{enumerate}
    Our next goal is to show that for all $i\in [C]$, there is some $P\in \calP'$ such that $P$ is mostly contained in $\mathbf{U}_i$. First note that since $|U_i|\geq 1$ for all $i\in [C]$, we have 
    \begin{align}\label{al:ci}
        n|[C]\setminus\{\sigma(P):P\in \calP'\}|&\leq\left|\bigcup_{i\in [C]\setminus \{\sigma(P):P\in \calP'\}}\mathbf{U}_i\right|.
    \end{align}
    Next, observe that 
    \begin{align}\label{al:vi}
  \left|\bigcup_{i\in [C]\setminus \{\sigma(P):P\in \calP'\}}\mathbf{U}_i\right|\leq  |V(\mathbf{H})|-\sum_{P\in \calP'}|P\cap \mathbf{U}_{\sigma(P)}|,
    \end{align}
where the inequality uses that the elements in $\calP'$ are pairwise disjoint.     We can use (i) and (ii) to obtain 
    \begin{align}\label{al:wi}
      \nonumber   \sum_{P\in \calP'}|P\cap \mathbf{U}_{\sigma(P)}| \geq  (1-\e^{1/4})\sum_{P\in \calP'} |P| &\geq (1-\e^{1/4})(1-2(k^2+1)\e^{1/2^k})|V(\mathbf{H})|\\
         &\geq (1-2(k^2+3)\e^{1/2^k})|V(\mathbf{H})|.
    \end{align}
    Finally, combining (\ref{al:ci}), (\ref{al:vi}), and (\ref{al:wi}),  we obtain
    $$
    n|[C]\setminus\{\sigma(P):P\in \calP'\}|\leq 2(k^2+3)\e^{1/2^k}|V(\mathbf{H})|.
    $$
    Dividing by $n$ yields
    $$
    |[C]\setminus\{\sigma(P):P\in \calP'\}|\leq 2(k^2+3)\e^{1/2^k}|V(H)|    < 1
    $$
where the last inequality is because  $|V(H)|\leq kC$ and $\e$ is sufficiently small. This implies $|\calP|\geq |\calP'|\geq C$.  This shows that there are arbitrarily large elements of $\calH$, whose smallest $\e$-regular partitions have at least $C$ parts.  Thus $M_{\calH}(\e)\geq C$.
\end{proof}

We now put everything together to prove our main theorem, Theorem \ref{thm:main}, restated below in more detail for the convenience of the reader. We have written the proof in such a way that one can easily deduce the claim preceding the statement of Theorem \ref{thm:main} concerning how each growth class is characterized. 

\begin{theorem}\label{thm:mainrestate}
Let $k\geq 3$ be an integer, and let $\calH$ be a hereditary $k$-graph property.  Then one of the following  holds.
\begin{enumerate}
\item (tower) There are $C_1,C_2>0$ so that for all sufficiently small $\e>0$,  
$$
Tw(\e^{-C_1})\leq M_{\calH}(\e)\leq Tw(\e^{-C_2}).
$$ 
\item (exponential) There are $C_1,C_2>0$ such that for all sufficiently small $\e>0$,  
$$
2^{\e^{-C_1}}\leq M_{\calH}(\e)\leq 2^{\e^{-C_2}};
$$ 
\item (polynomial) There are $C_1,C_2>0$ such that for all sufficiently small $\e>0$,  
$$
\e^{-C_1}\leq M_{\calH}(\e)\leq \e^{-C_2};
$$ 
\item (constant) There is $C>0$ such that for all sufficiently small $\e>0$,  
$$
M_{\calH}(\e)=C.
$$ 
\end{enumerate}
\end{theorem}
\begin{proof}
Before beginning, we set some notation. Given a class $\calF$ of $k$-graphs and an integer $N$, let $\calF[N]$ denote the class of all $N$-blowups of elements from $\calF$. Note the following holds by definition of ${\bf B}_{\calH}$, Fact \ref{fact:bhhgraphs}, and a simple pigeonhole argument.
\begin{align}\label{al:bh}
\text{For finite $\calF$, ${\bf B}_{\calH}\cap \calF\neq \emptyset$ if and only if $\calH\cap \calF[N]\neq \emptyset$ for all integers $N\geq 1$.}
\end{align}

We now prove the theorem. Suppose first that  $\calH\cap (n\otimes_k\calU_2(n))\neq \emptyset$ for all integers $n\geq 1$. Then by Theorem  \ref{thm:chung}, Proposition \ref{prop:svcequiv}, and Corollary \ref{cor:tower}(2) we have that for some constants $C_1,C_2>0$ and all sufficiently small $\e>0$,  
$$
Tw(\e^{-C_1})\leq M_{\calH}(\e)\leq Tw(\e^{-C_2}).
$$ 

Suppose now that for some integer $m\geq 1$, $\calH\cap (m\otimes_k\calU_2(m))= \emptyset$, but for all integers $n,N\geq 1$, $\calH\cap (\calU_k(n)[N])\neq  \emptyset$. By (\ref{al:bh}), for all integers $n\geq 1$, ${\bf B}_{\calH}\cap \calU_k(n)\neq \emptyset$, so by Proposition \ref{prop:vcclose}, $\calH$ is far from finite VC-dimension.  On the other hand, by Proposition \ref{prop:svcequiv}, $\calH$ is close to finite SVC-dimension.  By Proposition \ref{prop:polyexp}(2) and Corollary \ref{cor:tower}(1), there are constants $C_1,C_2>0$ such that for all sufficiently small $\e>0$,  
$$
2^{\e^{-C_1}}\leq M_{\calH}(\e)\leq 2^{\e^{-C_2}}.
$$ 
 Assume now that for all integers $n,N\geq 1$, $\calH\cap \calF_n^{(k)}[N]\neq \emptyset$, but for some integers $m,m'\geq 1$,  $\calH\cap (\calU_k(m)[m'])= \emptyset$. Then by (\ref{al:bh}), ${\bf B}_{\calH}\cap \calU_k(m)=\emptyset$, while ${\bf B}_{\calH}\cap \calF_n^{(k)}\neq \emptyset$ for all integers $n\geq 1$.  By Proposition \ref{prop:vcclose}, $\calH$ is close to finite VC-dimension. On the other hand, by definition of the $\calF_n^{(k)}$, the above implies ${\bf B}_{\calH}$ contains infinitely many non-isomorphic almost prime $k$-graphs. Combining with Proposition \ref{prop:polyexp}(1) and Proposition \ref{prop:polyconstant}(1), this implies that there are constants $C_1,C_2>0$ such that for all sufficiently small $\e>0$,  
$$
\e^{-C_1}\leq M_{\calH}(\e)\leq \e^{-C_2};
$$ 

Finally, we are left with the case where for some integers $n_1,n_2\geq 1$, $\calH\cap \calF_{n_1}^{(k)}[n_2]=\emptyset$.  By (\ref{al:bh}), this implies that for some $n_1$, ${\bf B}_{\calH}\cap \calF_{n_1}^{(k)}=\emptyset$.  Since ${\bf B}_{\calH}\neq \emptyset$ (see Fact \ref{fact:bhhgraphs} and recall we have implicitly assumed our hereditary property $\calH$ has arbitrarily large elements), we can conclude that in this remaining case, there is a unique maximal integer $K\geq 1$ such that ${\bf B}_{\calH}\cap \calF_K^{(k)}\neq \emptyset$.  By Theorem \ref{thm:constmain}, $M_{\calH}(\e)=K$ for all sufficiently small $\e>0$.  By (\ref{al:bh}) and Lemma \ref{lem:forb}, we see that this $K$ can be characterized as the unique integer satisfying that for all $N$, $\calH\cap \calF_K^{(k)}[N]\neq \emptyset$, but for some integer $m\geq 1$, $\calH\cap \calF_{K+1}^{(k)}[m]=\emptyset$.   
\end{proof}

\section{Appendix}

Here we prove Lemma \ref{lem:A(l)}, which we restate below for the convenience of the reader. We note the proof is essentially identical to the $3$-uniform case, proved in \cite{part2}, which in turn adapted well known arguments from model theory and structural graph theory.

\begin{lemma}\label{lem:A(l)repeat}
For all $\ell\geq 1$ there is $N\geq 1$ such that any almost prime $k$-graph with at least $N$ vertices contains an induced sub-$k$-graph which is an element of $\calA(\ell)$.
\end{lemma}
\begin{proof}
    Let $N$ be sufficiently large with respect to $k$ and $\ell$ and fix an almost prime $k$-graph $H$ with at least $N$ vertices. We will inductively construct an element of $\calA(\ell)$ as an induced sub-$k$-graph of $H$. We will use the following observation repeatedly:
    \begin{observation}\label{obs:not_sim}
        Since $H$ is almost prime, all $\sim_H$ classes have size at most $k$. Thus, if $X\subseteq V(H)$ has size at least $k+1$, there are $v_1,v_2\in X$ such that $v_1\not\sim_Hv_2$. 
    \end{observation}

    \underline{Step 1:} By Observation \ref{obs:not_sim}, there are $v_0, v_1\in V(H)$ such that $v_0\not\sim_Hv_1$. By definition of $\sim_H$ there are $x_1^2, \ldots, x_1^k\in V(H)\setminus\{v_0, v_1\}$ such that (after possibly relabeling) 
    \begin{align}\label{al:inout1}
    v_0x_1^2 \ldots x_1^k\notin E(H)\text{ and }v_1x_1^2 \ldots x_1^k\in E(H).
    \end{align}
     Define $\a(1)\in \{0,1\}$ to be such that 
    $$
    |N_{H^{\a(1)}}(x_1^2, \ldots, x_1^k)\setminus\{v_0, v_1\}|\geq \frac{|V(H)\setminus\{x_1^2, \ldots, x_1^k, v_0, v_1\}|}{2},
    $$
    and let $\beta(1)\in \{0,1\}\setminus\{\alpha(1)\}$. Define $z_1 = v_{\beta(1)}$ and 
    $$
    Z_1 = N_{H^{\a(1)}}(x_1^2, \ldots, x_1^k)\setminus\{v_0, v_1\}.
    $$
    We note the following.
    \begin{itemize}
        \item $\{\alpha(1), \beta(1)\} = \{0,1\}$
        \item $|\{z_1, x_1^2, \ldots, x_1^k\}| = k$ (using (\ref{al:inout1}))
        \item $z_1x_1^2\ldots x_1^k \in E^{\beta(1)}$
        \item $|Z_1|\geq \frac{|V(H)|-k-1}{2}\geq |V(H)|/4$.
    \end{itemize}
    This concludes Step 1.

    \underline{Step $t+1$}: Suppose we have defined the following:
    \begin{itemize}
        \item tuples $(z_1, x_1^2, \ldots, x_1^k), \ldots, (z_t, x_t^2, \ldots, x_t^k)$
        \item pairs $(\a(1), \beta(1)), \ldots, (\a(t), \beta(t))\in \{0,1\}^2$
        \item a set $Z_t\subset V(H)$
    \end{itemize}
    satisfying the following:
    \begin{enumerate}
        \item $\{\alpha(i), \beta(i)\} = \{0,1\}$ for all $i\in [t]$
        \item $|\{z_i, x_i^2, \ldots, x_i^k\}| = k$ for all $i\in [t]$
        \item $|\{z_1, \ldots, z_t\}|= t$ and $|\{(x_1^2, \ldots, x_1^k),\ldots,(x_t^2, \ldots, x_t^k)\}|= t$
        \item For all $i\in [t]$, $z_ix_i^2\ldots x_i^k\in E^{\beta(i)}$
        \item $Z_t$ satisfies  $|Z_t|\geq |V(H)|/4^{t}$ and
        $$
        Z_t\subseteq \left(\bigcap_{i=1}^tN_{H^{\a(i)}}(x_i^2,\ldots, x_i^k)\right)\setminus\{x_r^s, z_r:r\in[t], s\in \{2, \ldots, k\}\}.
        $$
       
    \end{enumerate}
    By Observation \ref{obs:not_sim}, there are $v_0, v_1\in Z_t$ such that $v_0\not\sim_H v_1$. Therefore, by definition of $\sim_H$, there is a $k-1$-tuple $(x_{t+1}^2, \ldots, x_{t+1}^k)$ such that (after possibly relabeling)  
    \begin{align}\label{al:inout2}
    v_0x_{t+1}^2\ldots x_{t+1}^k\notin E(H)\text{ and }v_1x_{t+1}^2\ldots x_{t+1}^k\in E(H).
    \end{align}
     Let $\a(t+1)$ be such that
    $$|N_{H^{\a(t+1)}}(x_{t+1}^2, \ldots, x_{t+1}^k)\cap Z_t\setminus\{v_0, v_1\}|\geq \frac{|Z_t\setminus\{x_{t+1}^2, \ldots, x_{t+1}^k, v_0, v_1\}|}{2}.$$
    Define $\beta(t+1)\in\{0,1\}\setminus\{\alpha(t+1)\}$, and let $z_{t+1} = v_{\beta(t+1)}$. Finally, define 
    $$Z_{t+1} = N_{H^{\alpha(t+1)}}(x_{t+1}^2, \ldots, x_{t+1}^k)\cap Z_t\setminus\{v_0, v_1\}$$
    We check that (1)-(5) still hold. 
    
    By how we have defined things, and by (\ref{al:inout2}), we see that (1), (2), and (4) hold for $t+1$.  The fact item (5) holds for $t$, since $v_0, v_1 \in Z_t$,  (\ref{al:inout2}) implies we must have $(x_{t+1}^2,\ldots, x_{t+1}^k)\neq(x_{i}^2,\ldots, x_{i}^k)$ for any $i\in[t]$. Therefore (3) holds for $t+1$. By construction and our induction hypotheses, we have 
    $$
    |Z_{t+1}|\geq \frac{|Z_t|-(k+2)}{2}\geq \frac{|Z_t|}{4}\geq \frac{|V(H)|}{4^{t+1}}.
    $$
     Further, since 
     $$
    Z_t\subseteq \left(\bigcap_{i=1}^tN_{H^{\a(i)}}(x_i^2,\ldots, x_i^k)\right)\setminus\{x_r^s, z_r:r\in[t], s\in \{2, \ldots, k\}\}
    $$
    and $Z_{t+1}\subseteq N_{H^{\alpha(t+1)}}(x_{t+1}^2, \ldots, x_{t+1}^k)\cap Z_t\setminus\{v_0, v_1\}$, we can conclude (5) holds for $t+1$. This completes Step $t+1$.

    Repeating this $s=k4^{2\ell+1}$-many times, which we may do since $N$ is sufficiently large, we have:
    \begin{itemize}
        \item A sequence of $(k-1)$-tuples of vertices $(x_1^2, \ldots, x_1^k), \ldots, (x_s^2, \ldots, x_s^k)$,
        \item A sequence of vertices $z_1, \ldots, z_s$, and
        \item A sequence $(\alpha(1), \beta(1)), \ldots, (\alpha(s), \beta(s))$,
    \end{itemize}
    such that
    \begin{itemize}
    \item for each $i\in [s]$, $\{\alpha(i),\beta(i)\}=\{0,1\}$,
        \item $(x_i^2, \ldots, x_i^k)\neq (x_j^2, \ldots, x_j^k)$ for each $i\neq j$ in $[s]$,
        \item $|\{z_i,x_i^2,\ldots, x_i^k\}|=k$ for each $i\in [s]$
        \item $z_i\neq z_j$ for all $i\neq j$ in $[s]$, 
        \item For all $i\in [s]$,
        $$
        z_ix_i^2\ldots x_i^k\in E^{\beta(i)}(H)
        $$
        and for all $1\leq i<j\leq s$, 
        $$
        z_jx_i^2\ldots x_i^k \in E^{\a(i)}(H).
        $$
    \end{itemize}

    By the pigeonhole principle, there is a choice of $\a\neq \beta\in \{0,1\}$, such that 
    $$
    |\{i\in [s]:(\a(i), \beta(i)) = (\a, \beta)\}|\geq \frac{s}{2}.
    $$
    We now inductively define a sequence of index sets. 
    
   At step $0$, we set $X_0 = \{i\in [s]:(\a(i), \beta(i)) = (\a, \beta)\}$. 
    
   At step $j\geq 1$, let $i_j = \max(X_{j-1})$. By pigeonhole principle there is $\tau(j)\in \{0,1\}$ such that 
    $$
    |N_{H^{\tau(j)}}(x_{i_j}^2, \ldots, x_{i_j}^k)\cap \{z_r:r\in X_{j-1}\setminus\{i_j\}\}|\geq \frac{|X_{j-1}|-k}{2}.
    $$
    We then set
    $$
    X_j = \{i\in X_{j-1}: i<i_j, z_i\in N_{H^{\tau(j)}}(x_{i_j}^2, \ldots, x_{i_j}^k)\},
    $$
    and note $|X_j|\geq |X_{j-1}|/4\geq |X_0|/4^{j}$.  
    
    Since $N$ is sufficiently large, we may repeat this process $2\ell$. From this we obtain indices $s\geq i_1> \ldots> i_{2\ell}\geq 1$ such that the following hold for all $j\in [2\ell]$:
    \begin{itemize}
        \item $z_{i_j}x_{i_j}^{2}\ldots x_{i_j}^k\in E^{\beta}$
        \item For all $1\leq r<j$, we have $z_{i_r}x_{i_j}^{2}\ldots x_{i_j}^k\in E^{\alpha}$
        \item For all $j<r\leq 2\ell$, we have $z_{i_r}x_{i_j}^{2}\ldots x_{i_j}^k\in E^{\tau(j)}$
    \end{itemize}
    Finally, by pigeonhole principle there is $\tau\in \{0,1\}$ such that $|\{j:\tau(j) = \tau\}|\geq \ell$. Let $X$ be an $\ell$-element subset of $\{j\in [2\ell]: \tau(j) = \tau\}$. Let  $j_1<\cdots<j_{\ell}$ be the elements of $X$, and for each $r\in [\ell]$, set $W_r=\{x_{i_{j_r}}^2,\ldots, x_{i_{j_r}}^k\}$. 
    
    If $(\alpha,\beta,\tau)\neq (1,0,0)$ then, after possibly relabeling, one can check the following induced sub-$k$-graph is in $\calA(\ell)$:
    $$
    H[\{z_{i_{j_1}}, \ldots,z_{i_{j_{\ell}}}\}\cup W_1\cup \ldots \cup W_{\ell}].
    $$
    If $(\alpha,\beta,\tau)=(1,0,0)$, let $W_1'=\{x_q^2,\ldots, x_q^k\}$ for any $q\in X_{2\ell}$ (which is nonempty by construction).  Then one can check the following induced sub-$k$-graph is in $\calA(\ell)$:
    $$
    H[\{z_{i_{j_1}}, \ldots,z_{i_{j_{\ell}}}\}\cup W_1'\cup W_2\cup \ldots \cup W_{\ell}].
    $$

\end{proof}

\end{document}